\documentclass[10pt, journal, twocolumn]{IEEEtran}

\pdfoutput=1

\ifCLASSINFOpdf
\else
\fi
\usepackage{float}
\usepackage{graphicx}
\usepackage{algorithm}
\usepackage{algpseudocode}
\usepackage{amsmath}
\usepackage{amssymb}
\usepackage{mathrsfs}
\usepackage{amsthm}
\usepackage{bm}
\usepackage{bbm}
\usepackage{cite}
\usepackage{enumitem}
\usepackage{diagbox}
\usepackage{color}
\usepackage{textgreek}
\usepackage{array}
\newcolumntype{L}[1]{>{\raggedright\let\newline\\\arraybackslash\hspace{0pt}}m{#1}}
\newcolumntype{C}[1]{>{\centering\let\newline\\\arraybackslash\hspace{0pt}}m{#1}}
\newcolumntype{R}[1]{>{\raggedleft\let\newline\\\arraybackslash\hspace{0pt}}m{#1}}
\usepackage[normalem]{ulem}  
\usepackage{verbatim}

\newcommand{\cK}{\mathcal{K}}
\newcommand{\cP}{\mathcal{P}}
\newcommand{\cN}{\mathcal{N}}

\newcommand{\cU}{\mathcal{U}}

\newcommand{\bw}{\boldsymbol w}
\newcommand{\bp}{\boldsymbol p}

\newcommand{\by}{\boldsymbol y}

\newcommand{\bx}{\boldsymbol x}

\newcommand{\bP}{\boldsymbol P}
\newcommand{\bzero}{\boldsymbol 0}

\newtheorem{myth}{Theorem}

\newtheorem{myle}[myth]{Lemma}
\newtheorem*{myre}{Remark}

\allowdisplaybreaks

\begin{document}
\title{Joint Subcarrier and Power Allocation in NOMA: Optimal and Approximate Algorithms} 

\author{Lou Sala\"un,~\IEEEmembership{Student Member,~IEEE}, Marceau~Coupechoux, and Chung~Shue~Chen,~\IEEEmembership{Senior Member,~IEEE}%
\thanks{L. Sala\"un and C. S. Chen are with Bell Labs, Nokia Paris-Saclay, 91620 Nozay, France, and also with Lincs, Paris 75013, France (e-mail: lou.salaun@nokia-bell-labs.com; chung\_shue.chen@nokia-bell-labs.com).}%
\thanks{L. Sala\"un and M. Coupechoux are with LTCI, Telecom ParisTech, University of Paris-Saclay, Paris 75013, France (e-mail: marceau.coupechoux@telecom-paristech.fr).}}

\maketitle

\begin{abstract}
Non-orthogonal multiple access (NOMA) is a promising technology to increase the spectral efficiency and enable massive connectivity in 5G and future wireless networks. In contrast to orthogonal schemes, such as OFDMA, NOMA multiplexes several users on the same frequency and time resource.
Joint subcarrier and power allocation problems (\textsc{JSPA}) in NOMA are NP-hard to solve in general.
In this family of problems, we consider the weighted sum-rate (WSR) objective function as it can achieve various tradeoffs between sum-rate performance and user fairness.
Because of JSPA's intractability, a common approach in the literature is to solve separately the power control and subcarrier allocation (also known as user selection) problems, therefore achieving sub-optimal result.  
In this work, we first improve the computational complexity of existing single-carrier power control and user selection schemes. These improved procedures are then used as basic building blocks to design new algorithms, namely \textsc{Opt-JSPA}, \textsc{\textepsilon-JSPA} and \textsc{Grad-JSPA}. \textsc{Opt-JSPA} computes an optimal solution with lower complexity than current optimal schemes in the literature. It can be used as a benchmark for optimal WSR performance in simulations. However, its pseudo-polynomial time complexity remains impractical for real-world systems with low latency requirements.
To further reduce the complexity, we propose a fully polynomial-time approximation scheme called \textsc{\textepsilon-JSPA}. Since, no approximation has been studied in the literature, \textsc{\textepsilon-JSPA} stands out by allowing to control a tight trade-off between performance guarantee and complexity.
Finally, \textsc{Grad-JSPA} is a heuristic based on gradient descent. Numerical results show that it achieves near-optimal WSR with much lower complexity than existing optimal methods.
\end{abstract}


\section{Introduction}\label{sec:intro}

In multi-carrier multiple access systems, the total frequency bandwidth is divided into subcarriers and assigned to users to optimize the spectrum utilization. Orthogonal multiple access (OMA), such as orthogonal frequency-division multiple access (OFDMA) adopted in 3GPP-LTE and also 5G New Radio Phase 1 standards~\cite{dahlman20185g,wong2017key}, only serves one user per subcarrier in order to avoid intra-cell interference and have low-complexity signal decoding at the receiver side. OMA is known to be suboptimal in terms of spectral efficiency~\cite{cover2012elements}.

The principle of multi-carrier non-orthogonal multiple access (MC-NOMA) is to multiplex several users on the same subcarrier by performing signals superposition at the transmitter side. Successive interference cancellation (SIC) is applied at the receiver side to mitigate interference between superposed signals. MC-NOMA is a promising multiple access technology for 5G and beyond mobile networks as it can achieve higher spectral efficiency than OMA schemes~\cite{saito2013system,dai2015non}. 

Careful optimization of the transmit powers is required to control the intra-carrier interference of superposed signals and maximize the achievable data rates. 
Besides, due to error propagation and decoding complexity concerns in practice~\cite{tse2005fundamentals}, subcarrier allocation for each transmission also needs to be optimized. 
As a consequence, joint subcarrier and power allocation problems (JSPA) in NOMA have received much attention. In this class of problems, weighted sum-rate (WSR) maximization is especially important as it can achieve different tradeoffs between sum-rate performance and user fairness~\cite{weeraddana2012weighted}.

Two types of power constraints are considered in the literature. On the one hand, cellular power constraint is mostly used in downlink transmissions to represent the total transmit power budget available at the base station (BS). On the other hand, individual power constraint sets a power limit independently for each user. The latter is often considered in uplink scenarios~\cite{chen2015suboptimal,al2014uplink}, nevertheless it can also be applied to the downlink~\cite{fu2018subcarrier,lei2016power}.

It is known that the equal-weight sum-rate and WSR maximization are both strongly NP-hard if we consider individual power constraints in OFDMA~\cite{liu2014complexity2} and in MC-NOMA systems~\cite{lei2016power,salaun2018optimal}. 
Nevertheless, several algorithms have been developed to perform subcarrier and/or power allocation for MC-NOMA and this type of constraints. Fractional transmit power control (FTPC) is a simple heuristic that allocates a fraction of the total power budget to each user based on their channel conditions~\cite{saito2013system}.
In~\cite{chen2015suboptimal} and~\cite{al2014uplink}, heuristic user pairing strategies and iterative resource allocation algorithms are studied for uplink transmissions.
A time efficient two-step heuristic is introduced in~\cite{fu2018subcarrier} to solve the problem with equal weights.
Reference~\cite{lei2016power} derives an upper bound on the optimal WSR and proposes a Lagrangian duality and dynamic programming (LDDP) scheme. This scheme achieves near-optimal result, assuming the power budget is divided in $J$ equal parts to be allocated. It mainly serves as benchmark due to its high computational complexity when $J$ is large in practical systems, which may not be suitable for low latency requirements. 

If we consider now cellular power constraints without individual power constraints, the equal-weight sum-rate maximization is polynomial time solvable in OFDMA~\cite{liu2014complexity1}. To the best of our knowledge, it is unknown whether WSR maximization in MC-NOMA under this type of constraints is polynomial time solvable or NP-hard. Reference~\cite[Proposition 1]{di2016subchannel} proves that the subcarrier optimization is NP-hard only in the case of equal power allocation among the users. The proposed polynomial-time reduction from the NP-complete 3-dimensional matching (3DM) to the NOMA problem should have shown that \underline{all} instances of 3DM can be mapped into an instance of the NOMA problem to be complete. Besides, the two-stage dynamic programming (TSDP) proposed in~\cite{lei2016power} solves it optimally in pseudo-polynomial time depending on $J$. Therefore, the WSR problem with cellular power constraint is \textit{weakly} NP-hard at most (in contrast to \textit{strongly} NP-hard for the individual power constraints as mentioned previously).
Only a few papers have developed optimization schemes in this setting, which are either heuristics with no theoretical performance guarantee or algorithms with impractical computational complexity.
For example, a greedy user selection and heuristic power allocation scheme based on difference-of-convex programming is proposed in~\cite{parida2014power}. In reference~\cite{di2016subchannel}, a matching algorithm is developed to perform subcarrier allocation. A minorization-maximization algorithm is used in~\cite{hanif2016minorization} to compute precoding vectors of a MISO-NOMA system. The authors of~\cite{sun2016optimal} employ monotonic optimization to develop an optimal resource allocation policy, which serves as benchmark due to its exponential complexity. The TSDP scheme is also optimal for cellular power constraint scenarios as proven in Theorem 13 of reference~\cite{lei2016power}, but it has high pseudo-polynomial complexity as well.

We note that, to the best of our knowledge, no polynomial-time approximation scheme (PTAS) has been proposed in the literature. Although PTAS is interesting for practical considerations of NP-hard problems, as it provides theoretical performance guarantees with controllable computational complexity. 
Motivated by this observation, we extend the framework of our previous paper~\cite{salaun2019weighted} with a fully polynomial-time approximation scheme (FPTAS) for the WSR maximization problem with cellular power constraint. In~\cite{salaun2019weighted}, we developed the following algorithms: two basic building blocks \textsc{SCPC} and \textsc{SCUS} which solve respectively the single-carrier power control and single-carrier user selection problems in polynomial time; and a heuristic JSPA scheme based on projected gradient descent, \textsc{SCPC} and \textsc{SCUS}, denoted here by \textsc{Grad-JSPA}. Our contributions are as follows:
\begin{enumerate}[leftmargin=*]
\item We improve \textsc{SCPC} and \textsc{SCUS} by performing precomputation to avoid repeated operations each time they are executed. This reduces their computational complexity by a factor proportional to the number of users.
\item The above precomputation also speeds up \textsc{Grad-JSPA}, which now has low and practical computational complexity. In addition, numerical results show that \textsc{Grad-JSPA} achieves near-optimal WSR, as well as significant improvement in performance over OMA.
\item We develop a new optimal algorithm, called \textsc{Opt-JSPA}, suitable for use as a benchmark for optimal WSR performance in simulations.
We show that \textsc{Opt-JSPA} has lower computational complexity than existing optimal schemes~\cite{lei2016power,sun2016optimal}.
\item We propose a FPTAS, which is denoted by \textsc{\textepsilon-JSPA}. Its design is based on the improved \textsc{SCPC} and \textsc{SCUS}, as well as techniques from the multiple choice knapsack problem~\cite{book:knapsack}. By definition of FPTAS, its performance is within a factor $1-\text{\textepsilon}$ of the optimal, for any $\text{\textepsilon}>0$. Moreover, it has polynomial complexity in both the input size and $1/\text{\textepsilon}$. Since, no approximation has been studied in the literature, \textsc{\textepsilon-JSPA} stands out by allowing to control a tight trade-off between performance guarantee and complexity.
\end{enumerate}

Through the aforementioned points, our aim is to deepen the understanding of \textsc{JSPA} and NOMA resource allocation problems. We develop optimal, approximate and heuristic schemes which are each suitable for systems with different computational capabilities, as well as for performance benchmarking.
In addition, we provide mathematical tools to study the WSR maximization problem, which can also be applied to other similar resource allocation problems.

The paper is organized as follows. In Section~\ref{sec:model}, we present the system model and notations. Section~\ref{sec:pbl} formulates the WSR problem. 
We consider two single-carrier sub-problems in Section~\ref{sec:SCopt} that were previously solved using \textsc{SCPC} and \textsc{SCUS} in~\cite{salaun2019weighted}. We propose improved versions of these algorithms, namely i-\textsc{SCPC} and i-\textsc{SCUS}, which perform precomputation to reduce their complexity. Based on these basic building blocks, we develop a low complexity gradient descent based heuristic (\textsc{Grad-JSPA}), a pseudo-polynomial time optimal algorithm (\textsc{Opt-JSPA}) and a FPTAS with \textepsilon-approximation guarantee (\textsc{\textepsilon-JSPA}) in Section~\ref{sec:MCopt}.
We show in Section~\ref{sec:num} some numerical results, highlighting our solution's WSR performance and computational complexity. In Section~\ref{sec:discussion}, we discuss about how to generalize our framework to more realistic channel estimation models and multi-antenna systems.
Finally, we conclude in Section~\ref{sec:ccl}.

\section{System Model and Notations}\label{sec:model}
We define in this section the system model and notations used throughout the paper. We consider a downlink multi-carrier NOMA system composed of one base station (BS) serving $K$ users. We denote the index set of users by $\cK \triangleq \{1,\ldots,K\}$, and the set of subcarriers by ${\cN \triangleq \{1,\ldots,N\}}$. The total system bandwidth $W$ is divided into $N$ subcarriers of bandwidth $W_n$, for each $n\in\cN$, such that $\sum_{n\in\cN}{W_n}=W$.
We assume orthogonal frequency division, so that adjacent subcarriers do not interfere each other. Moreover, each subcarrier $n\in\cN$ experiences frequency-flat block fading on its bandwidth $W_n$.

Let $p_k^n$ denotes the transmit power from the BS to user $k\in\cK$ on subcarrier $n\in\cN$. User $k$ is said to be active on subcarrier $n$ if $p_k^n > 0$, and inactive otherwise. In addition, let $g_k^n$ be the channel gain between the BS and user $k$ on subcarrier $n$, and $\eta_k^n$ be the received noise power. We assume that the channel gains are perfectly known. We discuss about more realistic models with imperfect channel state information (CSI) in Section~\ref{sec:discussion}.
For simplicity of notations, we define the normalized noise power as $\tilde{\eta}_k^n \triangleq \eta_k^n/g_k^n$. We denote by ${\bp \triangleq\left(p_k^n\right)_{k\in\cK,n\in\cN}}$ the vector of all transmit powers, and ${\bp^n \triangleq\left(p_k^n\right)_{k\in\cK}}$ the vector of transmit powers on subcarrier $n$.

In power domain NOMA, several users are multiplexed on the same subcarrier using superposition coding. A common approach adopted in the literature is to limit the number of superposed signals on each subcarrier to be no more than $M$. The value of $M$ is meant to characterize practical limitations of SIC due to decoding complexity and error propagation~\cite{tse2005fundamentals}.
We represent the set of active users on subcarrier $n$ by ${\cU_n \triangleq \{k\in\cK \colon p_k^n >0\}}$. The aforementioned constraint can then be formulated as $\forall n\in\cN,\; |\cU_n| \leq M$, where $|\mathord{\cdot}|$ denotes the cardinality of a finite set. Each subcarrier is modeled as a multi-user Gaussian broadcast channel~\cite{tse2005fundamentals} and SIC is applied at the receiver side to mitigate intra-band interference.

The SIC decoding order on subcarrier $n$ is usually defined as a permutation over the active users on $n$, i.e., $\pi_n \colon \{1,\ldots,|\cU_n|\} \to \cU_n$. However, for ease of reading, we choose to represent it by a permutation over all users $\cK$, i.e., $\pi_n \colon \{1,\ldots,K\} \to \cK$. These two definitions are equivalent in our model since the Shannon capacity~\eqref{eq:R_2} does not depend on the inactive users $k\in\cK\setminus\cU_n$, for which $p_k^n = 0$. For $i\in\{1,\ldots,K\}$, $\pi_n(i)$ returns the $i$-th decoded user's index. Conversely, user $k$'s decoding order is given by $\pi_n^{-1}(k)$. 

In this work, we consider the optimal decoding order studied in~\cite[Section 6.2]{tse2005fundamentals}. It consists of decoding users' signals from the highest to the lowest normalized noise power:
\begin{equation}\label{SIC_DL}
\tilde{\eta}_{\pi_n(1)}^n \geq \tilde{\eta}_{\pi_n(2)}^n \geq \cdots \geq \tilde{\eta}_{\pi_n(K)}^n.
\end{equation}
User $\pi_n(i)$ first decodes the signals of users $\pi_n(1)$ to $\pi_n(i-1)$ and subtracts them from the superposed signal before decoding its own signal. Interference from users $\pi_n(j)$ for $j > i$ is treated as noise. 
The maximum achievable data rate of user $k$ on subcarrier $n$ is given by Shannon capacity:
\begin{align} 
R_k^n(\bp^n) &\triangleq W_n\log_2\!\left(1+\frac{g_k^n p_k^n}{\sum_{j=\pi^{-1}_n(k)+1}^{K}{g_k^n p_{\pi_n(j)}^n}+\eta^n_k}\right), \nonumber\\
&\stackrel{\text{(a)}}{=} W_n\log_2\!\left(1+\frac{p_k^n}{\sum_{j=\pi^{-1}_n(k)+1}^{K}{p_{\pi_n(j)}^n}+\tilde{\eta}^n_k}\right), \label{eq:R_2}
\end{align}
where equality (a) is obtained after normalizing by $g_k^n$. We assume perfect SIC, therefore interference from users $\pi_n(j)$ for $j < \pi^{-1}_n(k)$ is completely removed in~\eqref{eq:R_2}. 

\section{Problem Formulation}\label{sec:pbl}

Let ${\bw = \{w_1,\ldots,w_K\}}$ be a sequence of $K$ positive weights. The main focus of this work is to solve the following JSPA optimization problem:

\begin{equation}\tag{$\cP$}\label{P}
\begin{aligned}
& \underset{\bp}{\text{maximize}}
& & \sum_{k\in\cK}{w_k\sum_{n\in\cN}{R_k^n\left(\bp^n\right)}}, \\
& \text{subject to}
& & C1:~\sum_{k\in\cK}\sum_{n\in\cN}{p_k^n} \leq P_{max}, \\
&&& C2:~\sum_{k\in\cK}{p_k^n} \leq P_{max}^n,~n\in\cN, \\
&&& C3:~p_k^n \geq 0,~k\in\cK,~n\in\cN, \\
&&& C4:~|\cU_n| \leq M,~n\in\cN. \\
\end{aligned}
\end{equation}
The objective of~\ref{P} is to maximize the system's WSR. As discussed in Section~\ref{sec:intro}, this objective function has received much attention since its weights $\bw$ can be chosen to achieve different tradeoffs between sum-rate performance and fairness~\cite{weeraddana2012weighted}. Note that $C1$ represents the cellular power constraint, i.e., a total power budget $P_{max}$ at the BS. In $C2$, we set a power limit of $P_{max}^n$ for each subcarrier $n$. This is a common assumption in multi-carrier systems, e.g.,~\cite{liu2014complexity1,liu2014complexity2}. Constraint $C3$ ensures that the allocated powers remain non-negative. Due to decoding complexity and error propagation in SIC~\cite{tse2005fundamentals}, we restrict the maximum number of multiplexed users per subcarrier to $M$ in $C4$. 

For ease of reading, we summarize some system parameters of a given instance of~\ref{P}, for all $n\in\cN$, as follows:
\begin{equation*}\label{param}
\mathcal{I}^n = \left(\bw,\cK,W_n,(g_k^n)_{k\in\cK},(\eta_k^n)_{k\in\cK}\right).
\end{equation*}

Let us consider the following change of variables:
\begin{equation}\label{eq:change}
\forall n\in\cN,\; x_i^n \triangleq \begin{cases}
    \sum_{j=i}^{K}{p_{\pi_n(j)}^n}, & \text{if $i\in\{1,\ldots,K\}$},\\
    0, & \text{if $i=K+1$}.
  \end{cases} 
\end{equation}
We define ${\bx \triangleq\left(x_i^n\right)_{i\in\{1,\ldots,K\},n\in\cN}}$ and ${\bx^n \triangleq\left(x_i^n\right)_{i\in\{1,\ldots,K\}}}$.
\begin{myle}[Equivalent problem~\ref{P2}]\label{le:equivalence}
$ $\newline
Problem~\ref{P} is equivalent to problem~\ref{P2} formulated below:
\begin{align}
& \underset{\bx}{\text{maximize}}
& & \sum_{n\in\cN}{\sum_{i=1}^{K}{f_i^n(x_i^n)}} + A, \tag{$\cP'$}\label{P2}\\
& \text{subject to}
& & C1':~\sum_{n\in\cN}{x_1^n} \leq P_{max}, \nonumber\\
&&& C2':~x_1^n \leq P_{max}^n,~n\in\cN, \nonumber\\
&&& C3':~x_i^n \geq x_{i+1}^n,~i\in\{1,\ldots,K\},~n\in\cN, \nonumber\\
&&& C3'':~x_{K+1}^n=0,~n\in\cN, \nonumber\\
&&& C4':~|\cU'_n| \leq M,~n\in\cN, \nonumber
\end{align}
where for any $i\in\{1,\ldots,K\}$ and $n\in\cN$, we have:
\begin{equation*}
f_i^n(x_i^n) \triangleq \begin{cases}
W_n\log_2\!\left(\left(x_1^n+\tilde{\eta}^n_{\pi_n(1)}\right)^{w_{\pi_n(1)}}\right), & \text{if $i=1$}, \\
W_n\log_2\!\left(\frac{\left(x_i^n+\tilde{\eta}^n_{\pi_n(i)}\right)^{w_{\pi_n(i)}}}{\left(x_i^n+\tilde{\eta}^n_{\pi_n(i-1)}\right)^{w_{\pi_n(i-1)}}}\right), & \text{if $i > 1$},
\end{cases}
\end{equation*}
and where ${\cU'_n \triangleq \{i\in\{1,\ldots,K\} \colon x_i^n > x_{i+1}^n\}}$. The constant term $A = \sum_{n\in\cN} w_{\pi_n(K)}\log_2\!\left(1/\tilde{\eta}^n_{\pi_n(K)}\right)$ is chosen so that~\ref{P} and~\ref{P2} have exactly the same optimal value.
\end{myle}
\begin{IEEEproof}
The idea is to apply the change of variables~\eqref{eq:change} to problem~\ref{P}. Details of the calculation can be found in Appendix~\ref{ap:transf}.
\end{IEEEproof}

The advantage of this formulation~\ref{P2} is that it exhibits a separable objective function in both dimensions $i\in\{1,\ldots,K\}$ and $n\in\cN$. In other words, it can be written as a sum of functions $f_i^n$, each only depending on one variable $x_i^n$. 

\section{Single-Carrier Optimization}\label{sec:SCopt}
In this section, we focus on a simpler problem, in which there is a single subcarrier $n\in\cN$ and a power budget $\bar{P}^n$ is given for this subcarrier:
\begin{align}
F^n\!\left(\bar{P}^n\right) = \; &\underset{\bx^n}{\text{max}}
& & {\sum_{i=1}^{K}{f_i^n\left(x_i^n\right)}} + A^n, \tag{$\cP'_{SC}(n)$}\label{Psub2}\\
& \text{subject to}
& & C2'\text{--}3':~\bar{P}^n \geq x_1^n \geq\ldots\geq x_K^n \geq 0,\nonumber\\
&&& C4':~|\cU'_n| \leq M, \nonumber
\end{align}
where $A^n = w_{\pi_n(K)}\log_2\!\left(1/\tilde{\eta}^n_{\pi_n(K)}\right)$. $C2'\text{--}3'$ is obtained by combining $C2'$, $C3'$ and $C3''$. $F^n\!\left(\bar{P}^n\right)$ denotes its optimal value.
Algorithms \textsc{SCPC} and \textsc{SCUS} have been introduced in our previous paper~\cite{salaun2019weighted} to tackle respectively the single-carrier power control and single-carrier user selection sub-problems that arise from~\ref{Psub2}. We provide technical details of these algorithms below, and we show how precomputation can further improve their computational complexity.
They will be used as basic building blocks in Section~\ref{sec:MCopt} to design efficient algorithms \textsc{Grad-JSPA}, \textsc{Opt-JSPA} and \textsc{\textepsilon-JSPA}, for the joint resource allocation problem. 

\subsection{Analysis of the Separable Functions $f_i^n$}\label{sec:ASF}
We introduce auxiliary functions to help us in the analysis of $f_i^n$ and the algorithm design. For $n\in\cN$, $i\in\{1,\ldots,K\}$ and $j \leq i$, assume that the consecutive variables $x_j^n,\ldots,x_i^n$ are all equal to a certain value $x\in\left[0,\bar{P}^n\right]$. We define $f_{j,i}^n$ as:
\begin{align}
f_{j,i}^n(x) &\triangleq \sum_{l=j}^{i}{f_l^n(x)}, \nonumber\\
&= \begin{cases}
W_n\log_2\!\left(\left(x+\tilde{\eta}^n_{\pi_n(i)}\right)^{w_{\pi_n(i)}}\right), & \text{if $j=1$}, \\
W_n\log_2\!\left(\frac{\left(x+\tilde{\eta}^n_{\pi_n(i)}\right)^{w_{\pi_n(i)}}}{\left(x+\tilde{\eta}^n_{\pi_n(j-1)}\right)^{w_{\pi_n(j-1)}}}\right), & \text{if $j > 1$}.\nonumber
\end{cases}
\end{align}
This simplification of notation is relevant for the analysis of \textsc{SCPC} and \textsc{SCUS} in the following subsections. Indeed, if users $j,\ldots,i-1$ are not active (i.e., $j,\ldots,i-1\notin\cU'_n$), then $x_j^n = \cdots = x_i^n$, therefore $\sum_{l=j}^{i}{f_l^n}$ can be replaced by $f_{j,i}^n$ and $x_{j+1}^n , \ldots , x_i^n$ are redundant with $x_j^n$. If constraint $C4'$ is satisfied, up to $M$ users are active on each subcarrier. Thus, evaluating the objective function of~\ref{Psub2} only requires $O\!\left(M\right)$ operations. 

\begin{figure}[!ht]
\vspace{-0.1cm}
\centering
\includegraphics[width=0.85\linewidth]{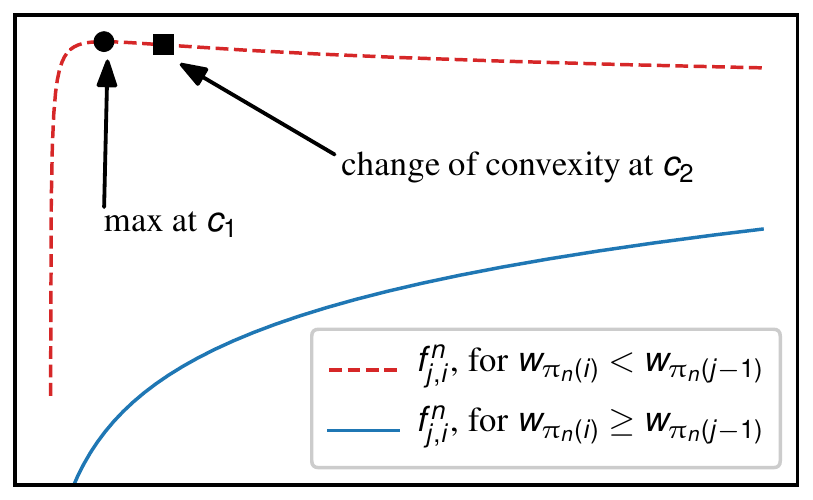}
\vspace{-0.2cm}
\caption{The two general forms of functions $f_{j,i}^n$}
\label{fig_f}
\vspace{-0.2cm}
\end{figure}

We study the properties of $f_{j,i}^n$ in Lemma~\ref{le:maxfij}. Note that $f_i^n = f_{i,i}^n$, therefore Lemma~\ref{le:maxfij} also holds for functions $f_i^n$. Fig.~\ref{fig_f} shows the two general forms that can be taken by $f_{j,i}^n$.

\begin{myle}[Properties of $f_{j,i}^n$]\label{le:maxfij}
$ $\newline
Let $n\in\cN$, $i\in\{1,\ldots,K\}$, and $j \leq i$, we have:
\begin{itemize}
\item If $j=1$ or $w_{\pi_n(i)} \geq w_{\pi_n(j-1)}$, then $f_{j,i}^n$ is increasing and concave on $\left[0,\infty\right)$.
\item Otherwise when $j>1$ and $w_{\pi_n(i)} < w_{\pi_n(j-1)}$, $f_{j,i}^n$ is unimodal. It increases on $\left(-\tilde{\eta}_{\pi_n(j-1)},c_1\right]$ and decreases on $\left[c_1,\infty\right)$, where
\begin{equation*}
c_1 = \frac{w_{\pi_n(j-1)}\tilde{\eta}_{\pi_n(i)}-w_{\pi_n(i)}\tilde{\eta}_{\pi_n(j-1)}}{w_{\pi_n(i)}-w_{\pi_n(j-1)}}.
\end{equation*}
Besides, $f_{j,i}^n$ is concave on $\left(-\tilde{\eta}_{\pi_n(j-1)},c_2\right]$ and convex on $\left[c_2,\infty\right)$, where
\begin{equation*}
c_2 = \frac{\sqrt{w_{\pi_n(j-1)}}\tilde{\eta}_{\pi_n(i)}-\sqrt{w_{\pi_n(i)}}\tilde{\eta}_{\pi_n(j-1)}}{\sqrt{w_{\pi_n(i)}}-\sqrt{w_{\pi_n(j-1)}}} \geq c_1.
\end{equation*}
\end{itemize}
\end{myle}
\begin{IEEEproof}
These analytical properties can be obtained by studying the first and second derivatives of $f_{j,i}^n$. Details can be found in Appendix~\ref{ap:propf}.
\end{IEEEproof}

\algnewcommand{\LineComment}[1]{\State \(\triangleright\) #1}
\newcommand{\algruleinput}[1][.2pt]{\par\vskip.2\baselineskip\hrule height #1\par\vskip.5\baselineskip}
\algnewcommand\algorithmicinput{\textbf{input:}}
\algnewcommand\Input{\item[\algorithmicinput]}
\algnewcommand\algorithmicglobalvar{\textbf{global variable:}}
\algnewcommand\GlobalVar{\item[\algorithmicglobalvar]}
\algnewcommand\algorithmicinit{\textbf{initialization:}}
\algnewcommand\Init{\item[\algorithmicinit]}
\algnewcommand\algorithmicfunc{\textbf{function}}
\algnewcommand\MyFunc{\item[\algorithmicfunc]}
\algnewcommand\algorithmicendfunc{\textbf{end function }}
\algnewcommand\EndFunc{\item[\algorithmicendfunc]}


\newlength{\textfloatsepsave} 
\setlength{\textfloatsepsave}{\textfloatsep} 
\setlength{\textfloatsep}{0.10cm}
\begin{algorithm}[ht]
\caption{Compute maximum of $f_{j,i}^n$ on $\left[0,\bar{P}^n\right]$}\label{algo:MaxF}
\begin{algorithmic}[1]
\MyFunc \textsc{Argmax}$f\!\left(j,i, \mathcal{I}^n, \bar{P}^n\right)$
\State $a \gets \pi_n(i)$
\State $b \gets \pi_n(j-1)$
\If{$j=1$ or $w_a \geq w_b$}
\State \Return $\bar{P}^n$
\Else
\State \Return $\max\left\{0,\min\left\{\frac{w_b\tilde{\eta}_a-w_a\tilde{\eta}_b}{w_a-w_b},\bar{P}^n\right\}\right\}$
\EndIf
\EndFunc
\end{algorithmic}
\end{algorithm}

We present in Algorithm~\ref{algo:MaxF} the pseudocode \textsc{Argmax}$f$ which computes the maximum of $f_{j,i}^n$ on $\left[0,\bar{P}^n\right]$ following the result of Lemma~\ref{le:maxfij}. \textsc{Argmax}$f$ only requires a constant number of basic operations, therefore its complexity is $O\!\left(1\right)$.

\subsection{Single-Carrier Power Control}\label{sec:SCPC}
The single-carrier power control problem~\ref{PsubSCPC} is equivalent to problem~\ref{Psub2}, with the exception that a fixed user selection $\cU'_n$ (or equivalently $\cU_n$) is given as input instead of being an optimization variable. It is defined below:
\begin{align}
&\underset{\bx^n}{\text{maximize}}
& & {\sum_{i=1}^{K}{f_i^n\left(x_i^n\right)}} + A^n, \tag{$\cP'_{SCPC}(n)$}\label{PsubSCPC1}\\
& \text{subject to}
& & C2'\text{--}3':~\bar{P}^n \geq x_1^n \geq\ldots\geq x_K^n \geq 0.\nonumber
\end{align}
We denote its optimal value by $F^n\!\left(\cU'_n,\bar{P}^n\right)$.

Since inactive users $k\notin\cU_n$ have no contribution on the data rates, i.e., $p_k^n = 0$ and $R_k^n = 0$, we remove them for the study of this sub-problem. Without loss of generality, we index the remaining active users on subcarrier $n$ by $i_n\in\{1_n,\ldots,|\cU'_n|_n\}$. For example, if $\cU'_n = \{4,7,10\}$, then $1_n = 4$, $2_n = 7$ and $3_n = 10$. For simplicity of notation, we add an index $0_n = 0$, which does not correspond to any user. From the definition of $\cU'_n$, variables $x_l^n$ with index from $l = (i-1)_n+1$ to $i_n$ are equal, for any $i \geq 1$. In the above example, we would have $x_1=x_2=x_3=x_4 > x_5 = x_6 = x_7 > x_8 = x_9 = x_{10}$. Thus, the objective function of~\ref{PsubSCPC1} can be written as:
\begin{equation}\label{activeUsersSimplification}
\sum_{i=1}^{K}{f_i^n(x_i^n)} + A^n = \sum_{i=1}^{|\cU'_n|} f_{(i-1)_n+1,i_n}^n\left(x_{i_n}^n\right) + B^n,
\end{equation}
where $B^n = A^n$ if the last active user's index is $|\cU'_n|_n=K$, and $B^n = f_{|\cU'_n|_n+1,K}^n\left(0\right) + A^n$ otherwise.
For $1 \leq j \leq i \leq K$, we simplify some notations as follows:
\begin{align*}
\tilde{f}_{j,i}^n\left(\cU'_n,\cdot\right) &\triangleq f_{(j-1)_n+1,i_n}^n\left(\cdot\right), \\
\textsc{Argmax}\tilde{f}\!\left(j,i,\mathcal{I}^n,\cU'_n,\bar{P}^n\right) &\triangleq\\
&\hspace{-1cm}\textsc{Argmax}f\!\left((j\!-\!1)_n\!+\!1,i_n, \mathcal{I}^n, \bar{P}^n\right).
\end{align*}
We reformulate the problem as:
\begin{align}
&\underset{x_{i_n}^n}{\text{maximize}} & & \sum_{i=1}^{|\cU'_n|}{\tilde{f}_{i,i}^n\left(\cU'_n,x_{i_n}^n\right)}+B^n,\tag{$\cP'_{SCPC}(n)$}\label{PsubSCPC}\\
& \text{subject to} & & C2'\text{--}3':~\bar{P}^n \geq x_{1_n}^n \geq\ldots\geq x_{|\cU'_n|_n}^n \geq 0.\nonumber
\end{align}
Algorithm~\ref{algoSCPC} presents the \textsc{SCPC} method proposed in~\cite{salaun2019weighted}.
The idea is to iterate over variables $x_{i_n}^n$ for $i = 1$ to $|\cU'_n|$, and compute their optimal value $x^* = \textsc{Argmax}\tilde{f}(i,i,\mathcal{I}^n,\cU'_n,\bar{P}^n)$ at line~3. If the current allocation satisfies constraint $C3'$, then $x_{i_n}^n$ gets value $x^*$. Otherwise, the algorithm backtracks at line~6 and finds the highest index $j \in\{1,\ldots,i-2\}$ such that $x_{j_n}^n \geq \textsc{Argmax}\tilde{f}(j+1,i,\mathcal{I}^n,\cU'_n,\bar{P}^n)$. Then, variables $x_{(j+1)_n}^n,\ldots,x_{i_n}^n$ are set equal to $\textsc{Argmax}\tilde{f}(j+1,i,\mathcal{I}^n,\cU'_n,\bar{P}^n)$ at line~10. The optimality and complexity of \textsc{SCPC} are presented in Theorem~\ref{th:algoSCPC}.

\begin{algorithm}[ht]
\caption{Single-carrier power control algorithm (\textsc{SCPC})}\label{algoSCPC}
\begin{algorithmic}[1]
\MyFunc \textsc{SCPC}$\left(\mathcal{I}^n, \cU'_n, \bar{P}^n\right)$
\For{$i=1$ to $|\cU'_n|$}
\LineComment Compute the optimal of $\tilde{f}_{i,i}^n$
\State $x^* \gets \textsc{Argmax}\tilde{f}\!\left(i,i,\mathcal{I}^n,\cU'_n,\bar{P}^n\right)$
\LineComment Modify $x^*$ if this allocation violates constraint $C3'$
\State $j \gets i-1$
\While{$x_{j_n}^n < x^*$ and $j\geq 1$}
\State $x^* \gets \textsc{Argmax}\tilde{f}\!\left(j,i,\mathcal{I}^n,\cU'_n,\bar{P}^n\right)$
\State $j \gets j-1$
\EndWhile
\State $x_{(j+1)_n}^n,\ldots,x_{i_n}^n \gets x^*$
\EndFor
\State \Return $x_{1_n}^n,\ldots,x_{|\cU'_n|_n}^n$
\EndFunc
\end{algorithmic}
\end{algorithm}

\begin{myth}[Optimality and complexity of \textsc{SCPC}]\label{th:algoSCPC}
$ $\newline
Given subcarrier $n\in\cN$, a set $\cU'_n$ of $M$ active users and a power budget $\bar{P}^n$, algorithm \textsc{SCPC} computes the optimal single-carrier power control. Its worst case computational complexity is $O\!\left(M^2\right)$. 
\end{myth}
\begin{IEEEproof}
We prove this theorem in Appendix~\ref{ap:thalgoSCPC} by mathematical induction combined with Lemma~\ref{le:maxfij}.
\end{IEEEproof}

In multi-carrier resource allocation schemes, such as \textsc{Grad-JSPA} and \textsc{\textepsilon-JSPA}, it is often required to compute the optimal single-carrier power control and WSR for many different values of power budget $\bar{P}^n$.
In these cases, running many times \textsc{SCPC} is actually not efficient in terms of computational complexity, since several computations may be repeated. To avoid this, we propose in Algorithm~\ref{algoSCPC_speed} an improved \textsc{SCPC} algorithm (i-\textsc{SCPC}). The idea is to perform precomputation before runtime by calling \textsc{SCPC}$\left(\mathcal{I}^n, \cU'_n, P_{max}\right)$ and storing its result $x_{1_n}^n,\ldots,x_{|\cU'_n|_n}^n$ as a global variable (also called lookup table). 
Any subsequent evaluation with input $\mathcal{I}^n$, $\cU'_n$, $\bar{P}^n$, where $\bar{P}^n\leq P_{max}$, can be obtained as in line~1.

\begin{algorithm}[ht]
\caption{Improved \textsc{SCPC} algorithm with precomputation}\label{algoSCPC_speed}
\begin{algorithmic}[1]
\Input $\mathcal{I}^n, \cU'_n, P_{max}$
\GlobalVar $x_{1_n}^n,\ldots,x_{|\cU'_n|_n}^n$
\Init $x_{1_n}^n,\ldots,x_{|\cU'_n|_n}^n \gets$  \textsc{SCPC}$\left(\mathcal{I}^n, \cU'_n, P_{max}\right)$
\MyFunc i-\textsc{SCPC}$\left(\bar{P}^n\right)$
\State \Return $\min\{x_{1_n}^n,\bar{P}^n\},\ldots,\min\{x_{|\cU'_n|_n}^n,\bar{P}^n\}$
\EndFunc
\end{algorithmic}
\end{algorithm}

\begin{myth}[Optimality and complexity of i-\textsc{SCPC}]\label{th:algoiSCPC}
$ $\newline
Given subcarrier $n\in\cN$ and a set $\cU'_n$ of $M$ active users, the precomputation of i-\textsc{SCPC} has complexity $O\!\left(M^2\right)$. Any subsequent evaluation costs $O\!\left(M\right)$. Hence, for $C$ different power budgets, algorithm i-\textsc{SCPC} computes their respective optimal single-carrier power control with overall complexity $O\!\left(M^2+CM\right)$.
\end{myth}
\begin{IEEEproof}
We show in Appendix~\ref{ap:algoiSCPC} that subsequent evaluations of \textsc{SCPC} can be obtained as in line~1 of Algorithm~\ref{algoSCPC_speed}.
\end{IEEEproof}

\begin{myre}
Note that \textsc{SCPC} and i-\textsc{SCPC} returns $|\cU'_n|$ values $x_{1_n}^n, \ldots, x_{|\cU'_n|_n}^n$ representing only the active users' variables. These values are sufficient to compute the optimal power allocation and WSR of~\ref{PsubSCPC} as shown in Eqn.~(\ref{activeUsersSimplification}). If needed, the full vector $\bx^n$ can be obtained by the following procedure in $O\!\left(K\right)$ operations: 
\begin{algorithmic}[1]
\For{$i = 1$ to $|\cU'_n|$ and $l = (i-1)_n+1$ to $i_n$}
\State $x_l^n \gets x_{i_n}^n$
\EndFor
\For{$l = |\cU'_n|_n+1$ to $K$}
\State $x_l^n \gets 0$
\EndFor
\end{algorithmic}

\end{myre}

\subsection{Single-Carrier User Selection}\label{sec:SCUS}

Unlike in the previous subsection, we consider here furthermore the user selection $\cU'_n$ optimization under multiplexing and SIC constraint $C4'$, i.e., we solve  \ref{Psub2}.
In~\cite{salaun2019weighted}, a dynamic programming (DP) is proposed to solve~\ref{Psub2}. Here, we first develop a similar DP procedure in Algorithm~\ref{alg3} (\textsc{SCUS}). Then, we propose an improved version (i-\textsc{SCUS}) which performs \textsc{SCUS} as precomputation. 

The idea of \textsc{SCUS} is to compute recursively the elements of three arrays $V$, $X$, $U$. Let $m\in\{0,\ldots,M\},\;j\in\{1,\ldots,K\}$ and $i\in\{j,\ldots,K\}$, we define $V[m,j,i]$ as the optimal value of the following problem~\ref{defV}:
\begin{align}
V[m,j,i] & & \triangleq &&& \underset{\bx^n}{\max} & & {\sum_{l=j}^{K}{f_l^n(x_l^n)}}, \tag{$\cP'_{SC}[m,j,i]$}\label{defV}\\
&&&&& \text{subject to} & & C2', C3', C3'',\nonumber\\
&&&&&&& C4':~|\cU'_n| \leq m, \nonumber\\
&&&&&&& C5':~x_j^n = \cdots = x_i^n. \nonumber
\end{align}
This problem is more restrictive than~\ref{Psub2}. The objective function only depends on variables $x_j^n,\ldots,x_K^n$. $C4'$ limits the number of active users to $m$. Moreover, variables $x_j^n,\ldots,x_i^n$ are equal according to $C5'$.

\begin{algorithm}[!t]
\caption{Single-carrier user selection algorithm (\textsc{SCUS})}\label{alg3}
\begin{algorithmic}[1]
\MyFunc \textsc{SCUS}$\left(\mathcal{I}^n, M, \bar{P}^n\right)$
\LineComment Initialize arrays $V$, $X$, $U$ for $m=0$ and $i=K$
\For{$i = K$ to $1$ and $j = i$ to $1$}
\State $V\left[0,j,i\right] \gets f_{j,K}^n\left(0\right)$
\State $X\left[0,j,i\right] \gets 0$
\State $U\left[0,j,i\right] \gets \varnothing$
\EndFor
\For{$m = 1$ to $M$ and $j = K$ to $1$}
\State $x^* \gets \textsc{Argmax}f\!\left(j,K,\mathcal{I}^n,\bar{P}^n\right)$
\State $V\left[m,j,K\right] \gets f_{j,K}^n\left(x^*\right)$
\State $X\left[m,j,K\right] \gets x^*$
\State $U\left[m,j,K\right] \gets \varnothing$
\EndFor
\LineComment Compute $V$, $X$, $U$ for $m\in\left[1,M\right]$ and $j\leq i \leq K-1$
\For{$i = K-1$ to $1$ and $m = 1$ to $M$ and $j = i$ to $1$}
\State $x^* \gets \textsc{Argmax}f\!\left(j,i,\mathcal{I}^n,\bar{P}^n\right)$
\State $v_{act} \gets f_{j,i}^n\left(x^*\right) + V[m-1,i+1,i+1]$
\State $v_{inact} \gets V[m,j,i+1]$
\If{$v_{act} > v_{inact}$ and $x^* > X\left[m-1,i+1,i+1\right]$}
\State $V\left[m,j,i\right] \gets v_{act}$
\State $X\left[m,j,i\right] \gets x^*$
\State $U\left[m,j,i\right] \gets (m-1,i+1,i+1)$
\Else
\State $V\left[m,j,i\right] \gets v_{inact}$
\State $X\left[m,j,i\right] \gets  X\left[m,j,i+1\right]$
\State $U\left[m,j,i\right] \gets (m,j,i+1)$
\EndIf
\EndFor
\LineComment Retrieve the optimal solution $\bx^n$
\State $x_1^n,\ldots,x_K^n \gets 0$
\State $(m,j,i) \gets (M,1,1)$
\Repeat 
\State $x_j^n,\ldots,x_i^n \gets X\left[m,j,i\right]$
\State $(m,j,i) \gets U\left[m,j,i\right]$
\Until{$(m,j,i) = \varnothing$}
\State\Return $\bx^n$
\EndFunc
\end{algorithmic}
\end{algorithm}

It is interesting to note that $V[M,1,1]$ is the optimal value of~\ref{Psub2}, since the objective function is $\sum_{l=1}^{K}{f_l^n(x_l^n)}$ for $j=1$ and constraint $C5'$ becomes trivially true for $j = i$. Let ${x_j^n}^*,\ldots,{x_K^n}^*$ be the optimal solution achieving $V[m,j,i]$. We define $X[m,j,i] \triangleq {x_i^n}^*$, which is also equal to ${x_j^n}^*,\ldots,{x_{i-1}^n}^*$ due to constraint $C5'$. The idea of SCUS is to recursively compute the elements of $V$ through the following relation:
\begin{equation}\label{SCUSrec}
V[m,j,i] = \begin{cases}
v_{act}, \quad \text{if $v_{act} > v_{inact}$} \\
\hphantom{v_{inact},} \text{ and $x^* > X\left[m-1,i+1,i+1\right]$}, \\
v_{inact}, \text{ otherwise},
\end{cases}
\end{equation}
where $x^* = \textsc{Argmax}f\!\left(j,i, \mathcal{I}^n, \bar{P}^n\right)$, and $v_{act}$ (resp. $v_{inact}$) corresponds to allocation where user $i$ is active (resp. inactive):
\begin{align*}
&v_{act} = f_{j,i}^n\left(x^*\right) + V[m-1,i+1,i+1], \\
&v_{inact} = V[m,j,i+1].
\end{align*}

During SCUS's iterations, the array $U$ keeps track of which previous element of $V$ has been used to compute the current value function $V[m,j,i]$. This allows us to retrieve the entire optimal vector $\bx^n$ at the end of Algorithm~\ref{alg3} (at lines 28-35) by backtracking from index $(M,1,1)$ to $\varnothing$, where $\varnothing$ is set at initial indices (see lines 5 and 11) to indicate the recursion termination. To sum up, $X$ and $U$ have two different recurrence relations depending on the cases in Eqn.~(\ref{SCUSrec}).\\
If $V\left[m,j,i\right] = v_{act}$, then:
\begin{align*}
X\left[m,j,i\right] &= x^*,\\
U\left[m,j,i\right] &= (m-1,i+1,i+1).
\end{align*}
If $V\left[m,j,i\right] = v_{inact}$, then:
\begin{align*}
X\left[m,j,i\right] &= X\left[m,j,i+1\right],\\
U\left[m,j,i\right] &= (m,j,i+1).
\end{align*}

When $m = 0$, no user can be active on this subcarrier due to constraint $C4'$. Therefore, $V$, $X$, $U$ can be initialized by: 
\begin{align*}
V\left[0,j,i\right] &= f_{j,K}^n\left(0\right),\\
X\left[0,j,i\right] &= 0,\\
U\left[0,j,i\right] &= \varnothing.
\end{align*}
For simplicity, we also extend $V$, $X$ and $U$ on the index $i = K$ and $j \leq K$ and initialize them as follows:
\begin{align*}
V\left[m,j,K\right] &= f_{j,K}^n\left(x^*\right),\\
X\left[m,j,K\right] &= x^*,\\
U\left[m,j,K\right] &= \varnothing.
\end{align*}
A detailed analysis is given in Appendix~\ref{ap:thalgoSCUS}.

\begin{myth}[Optimality and complexity of \textsc{SCUS}]\label{th:algoSCUS}
$ $\newline
Given a subcarrier $n\in\cN$, a power budget $\bar{P}^n$ and $M\geq 1$, algorithm \textsc{SCUS} computes the optimal single-carrier power control and user selection of~\ref{Psub2}. Its worst case computational complexity is $O\!\left(MK^2\right)$.
\end{myth}
\begin{IEEEproof}
The proof is done by induction based on the principle of dynamic programming. See Appendix~\ref{ap:thalgoSCUS}.
\end{IEEEproof}

We present i-\textsc{SCUS} in Algorithm~\ref{algiSCUS}, which performs precomputation to avoid repeating the DP procedure when multiple evaluations are required. The algorithm precomputes vectors $V$, $X$, $U$ from \textsc{SCUS}$\left(\mathcal{I}^n, M, P_{max}\right)$ before runtime, at line~1. Then, in lines~2-5, it retrieves the active users set $\cU'_n$ and optimal solution $x_{1}^n,\ldots,x_{K}^n$ of each $V[M,1,i]$, $i\in\{1,\ldots,K\}$, and stores them in $collection$. Any subsequent evaluation with a lower budget $\bar{P}^n\leq P_{max}$, can be obtained by searching the best allocation among the $K$ possibilities in $collection$ (lines~6-7). Each allocation is truncated as in i-\textsc{SCPC}$\left(\bar{P}^n\right)$ to satisfy budget $\bar{P}^n$. The optimality and complexity of Algorithm~\ref{algiSCUS} are given in Theorem~\ref{th:algoiSCUS}.

\begin{algorithm}[ht]
\caption{Improved \textsc{SCUS} algorithm with precomputation}\label{algiSCUS}
\begin{algorithmic}[1]
\Input $\mathcal{I}^n, M, P_{max}$
\GlobalVar $collection$
\Init 
\State Get $V,X,U$ from \textsc{SCUS}$\left(\mathcal{I}^n, M, P_{max}\right)$
\For{$i=1$ to $K$}
\State Retrieve the active users set $\cU'_n$ of $V[M,1,i]$ and its $\hphantom{for}$optimal solution $x_{1}^n,\ldots,x_{K}^n$
\State Add $\left(\cU'_n, x_{1}^n,\ldots,x_{K}^n\right)$ to $collection$
\EndFor
\MyFunc i-\textsc{SCUS}$\left(\bar{P}^n\right)$
\State Get $\left(\cU'_n, x_{1}^n,\ldots,x_{K}^n\right)$ in $collection$ that maximizes $F^n(\cU'_n,\bar{P}^n) = \sum_{l=1}^{|\cU'_n|} \tilde{f}_{l,l}^n\left(\cU'_n,\min\{x_{l_n}^n,\bar{P}^n\}\right) + B^n$
\State \Return $ \min\{x_{1}^n,\bar{P}^n\},\ldots,\min\{x_{K}^n,\bar{P}^n\}$
\EndFunc
\end{algorithmic}
\end{algorithm}

\begin{myth}[Optimality and complexity of i-\textsc{SCUS}]\label{th:algoiSCUS}
$ $\newline
Given a subcarrier $n\in\cN$, a power budget $\bar{P}^n$ and $M\geq 1$, the precomputation of i-\textsc{SCUS} has complexity $O\!\left(MK^2\right)$. Any subsequent evaluation costs $O\!\left(MK\right)$. Hence, for $C$ different power budgets, i-\textsc{SCUS} computes their respective optimal single-carrier power control and user selection~\ref{Psub2} with overall complexity $O\!\left(MK^2+CMK\right)$.
\end{myth}
\begin{IEEEproof}
See Appendix~\ref{ap:algoiSCUS}.
\end{IEEEproof}

Table~\ref{table:SC} summarizes the complexity of the single-carrier algorithms developed in this section. They will be used as basic building blocks to design JSPA schemes in Section~\ref{sec:MCopt}.
\begin{table}[ht]
\begin{center}
\caption{Summary of the single-carrier resource allocation schemes}
\vspace{-0.1cm}
\label{table:SC}
\begin{tabular}{|c|c|}
  \hline
  \textbf{Algorithm} & \textbf{Complexity to perform $C$ evaluations} \\ \hline
  \textsc{SCPC}~\cite{salaun2019weighted} & $O\!\left(CM^2\right)$\\\hline
  i-\textsc{SCPC} & $O\!\left(M^2+CM\right)$\\\hline
  \textsc{SCUS}~\cite{salaun2019weighted} & $O\!\left(CMK^2\right)$\\\hline
  i-\textsc{SCUS} & $O\!\left(MK^2+CMK\right)$\\\hline
\end{tabular}
\end{center}
\end{table}

\section{Joint Subcarrier and Power Allocation}\label{sec:MCopt}
Recall that $F^n\!\left(\bar{P}^n\right)$ is the optimal value of~\ref{Psub2} with power budget $\bar{P}^n$. We have $F^n\!\left(\bar{P}^n\right) = \sum_{i=1}^{K}{f_i^n(x_i^n)} + A^n$, where $x_{1}^n,\ldots,x_{K}^n$ is the output of i-\textsc{SCUS}$\left(\bar{P}^n\right)$.
Using this notation, the JSPA problem~\ref{P2} can be simplified as:
\begin{align}
& \underset{\bar{\bP}}{\text{maximize}}
& & \sum_{n\in\cN}F^n\!\left(\bar{P}^n\right), \tag{$\cP'_{MC}$}\label{Psub1}\\
& \text{subject to}
& & \bar{P}^n \in \mathcal{F}, \nonumber
\end{align}
where $\bar{P}^n$, for $n\in\cN$, are intermediate variables representing each subcarrier's power budget. $\bar{\bP} \triangleq \left(\bar{P}^1,\ldots,\bar{P}^N\right)$ denotes the power budget vector. The feasible set 
\begin{equation*}
\mathcal{F} \triangleq \{\bar{\bP} \colon \sum_{n\in\cN}\bar{P}^n \leq P_{max} \text{ and }  0 \leq \bar{P}^n \leq P_{max}^n,\; n\in\cN\}
\end{equation*} 
is chosen to satisfy $C1'$ and $C2'$ in~\ref{P2}.

Problem~\ref{Psub1} consists in optimizing the power budget $\bar{P}^n$ allocated to each subcarrier $n$. For a given budget $\bar{P}^n$, $F^n\!\left(\bar{P}^n\right)$ is computed by finding the optimal single-carrier user selection and power control using i-\textsc{SCUS}$\left(\bar{P}^n\right)$. The choice of $\bar{P}^n$ affects the single-carrier user selection and power control, i.e., variables $\bx^n$ and $\cU'_n$. The latter influence the value of $F^n$, which in turn has an impact on the power budget optimization. Although variables $\bx^n$ and $\cU'_n$ are hidden in $F^n$, they are nevertheless optimized jointly with $\bar{P}^n$. Indeed, we can see that~\ref{Psub1} is equivalent to~\ref{P2} when replacing $F^n$ by its definition in~\ref{Psub2} along with its constraints $C2'$ to $C4'$.

\vspace{-0.3cm}

\subsection{Gradient Descent Based Heuristic}\label{sec:gradJSPA}
\begin{figure}[!ht]
\vspace{-0.1cm}
\centering
\includegraphics[trim={0 0.5cm 0 0.5cm},width=0.95\linewidth]{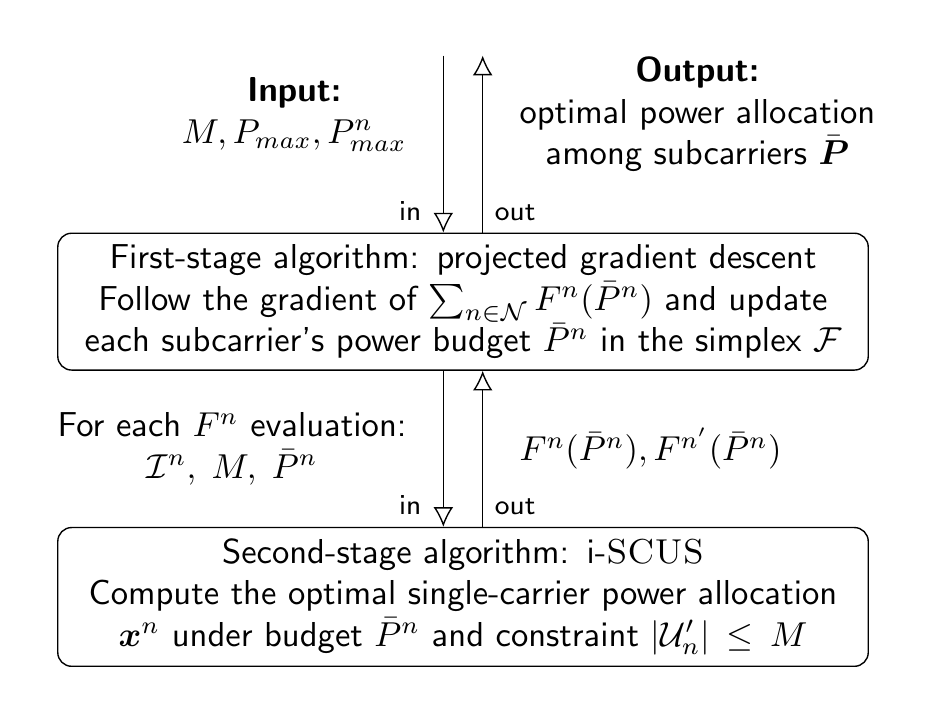}
\caption{Overview of \textsc{Grad-JSPA}}
\label{fig_diagramGradJSPA}
\end{figure}

\textsc{Grad-JSPA} is an efficient heuristic based on projected gradient descent. Its principle is to perform a two-stage optimization as presented in Fig.~\ref{fig_diagramGradJSPA}. The first-stage is a projected gradient descent on $\bar{\bP}$ in the search space $\mathcal{F}$. The gradient descent requires to evaluate $\sum_{n\in\cN}F^n$ and its gradient at each iteration. This task is carried out by i-\textsc{SCUS} 
in the second-stage.
We denote the derivative of $F^n$ at $\bar{P}^n$ by $F^{n\prime}\!\left(\bar{P}^n\right)$. Lemma~\ref{le:derivative} shows how to compute it. As illustrated in Fig.~\ref{fig_diagramGradJSPA}, the second-stage is called at each gradient iteration to return $F^{n}\!\left(\bar{P}^n\right)$ and $F^{n\prime}\!\left(\bar{P}^n\right)$ to the first-stage, for all $n\in\cN$.

\begin{myle}[Derivative of $F^n$]\label{le:derivative}
$ $\newline
Let $x_{1}^n,\ldots,x_{K}^n$ be the output of \mbox{i-\textsc{SCUS}$\left(\bar{P}^n\right)$}. The left derivative of $F^n$ at $\bar{P}^n$, can be computed as follows:
\begin{equation*}
F^{n\prime}\!\left(\bar{P}^n\right) = \frac{W_n w_{\pi_n(l)}}{\left(x_l^n+\tilde{\eta}^n_{\pi_n(l)}\right)\ln\!{(2)}} = \frac{W_n w_{\pi_n(l)}}{\left(\bar{P}^n+\tilde{\eta}^n_{\pi_n(l)}\right)\ln\!{(2)}},
\end{equation*}
where $l$ is the greatest index such that $x_l^n = \bar{P}^n$, and $\ln\!{(2)}$ is the natural logarithm of $2$.
\end{myle}
\begin{IEEEproof}
To get $F^{n\prime}$, we first prove that $F^{n\prime}\!\left(\cU'_n,\bar{P}^n\right) \!=\! f_{1,l}^{n\prime}\!\left(\bar{P}^n\right)$
 and $F^n(\bar{P}^n) \!=\! \max_{\,\cU'_n}\{F^n(\cU'_n,\bar{P}^n)\}$, where $\max$ is taken over all active users sets in $collection$ of \mbox{i-\textsc{SCUS}}.
See Appendix~\ref{ap:derivative} for the detailed proof of semi-differentiability.
\end{IEEEproof}

The pseudocode of \textsc{Grad-JSPA} is described in Algorithm~\ref{algoGradJSPA}. 
Input $\xi$ corresponds to the error tolerance at termination, as we can see at line~8. 
The search direction at lines~4-5 is the gradient of $\sum_{n\in\cN}F^n$ evaluated at $\bar{\bP}$. Since $F^1,\ldots,F^N$ are independent, it is equal to the vector of $F^{1\prime}\!\left(\bar{P}^1\right),\ldots,F^{N\prime}\!\left(\bar{P}^N\right)$.
Note that the step size $\alpha$ at line~6 can be tuned by backtracking line search or exact line search~\cite[Section 9.2]{boyd2004convex}. We adopt the latter to perform simulations. The projection of $\bar{\bP} + \alpha\Delta$ on the simplex $\mathcal{F}$ at line~7 can be computed efficiently~\cite[Section 8.1.1]{boyd2004convex}, the details of its implementation are omitted here.

\begin{algorithm}[t]
\caption{Gradient descent based heuristic (\textsc{Grad-JSPA})}\label{algoGradJSPA}
\begin{algorithmic}[1]
\MyFunc \textsc{Grad-JSPA}$\left(\left(\mathcal{I}^n\right)_{n\in\cN},M,P_{max},P^n_{max},\xi\right)$
\State Let $\bar{\bP} \gets \bzero$ be the starting point
\Repeat
\State Save the previous vector $\bar{\bP}' \gets \bar{\bP}$
\State Determine a search direction $\Delta$
\State \hspace*{5mm} $\Delta \gets \left(F^{1\prime}\!\left(\bar{P}^1\right),\ldots,F^{N\prime}\!\left(\bar{P}^N\right)\right)$
\State Choose a step size $\alpha$
\State Update $\bar{\bP} \gets$ projection of $\bar{\bP} + \alpha\Delta$ on $\mathcal{F}$
\Until{$||\bar{\bP}'-\bar{\bP}||\leq \xi$}
\State \Return $\bar{\bP}$
\EndFunc
\end{algorithmic}
\end{algorithm}

We showed in our previous work~\cite{salaun2019weighted} that \textsc{Grad-JSPA} worst case complexity is $O\!\left(\log\!\left(1/\xi\right)\!NMK^2\right)$ when \textsc{SCUS} is used to evaluate functions $F^n$, $n\in\cN$. We show in Theorem~\ref{th:algoGradJSPA} that the complexity of \textsc{Grad-JSPA} can be reduced by the use of i-\textsc{SCUS}.

\begin{myth}[Complexity of \textsc{Grad-JSPA}]\label{th:algoGradJSPA}
$ $\newline
Let $\xi$ be the error tolerance at termination. Algorithm \textsc{Grad-JSPA} has complexity $O\!\left(NMK^2+\log\!\left(1/\xi\right)\!NMK\right)$ when i-\textsc{SCUS} is used to evaluate functions $F^n$, $n\in\cN$. 
\end{myth}
\begin{IEEEproof}
In Appendix~\ref{ap:derivative}, we prove that the objective function is piece-wise $\alpha$-strongly concave and $\beta$-smooth. Therefore, the convergence to a local optimum follows from classical convex optimization results~\cite[Section 2.2.4]{nesterov1998introductory}.
\end{IEEEproof}

Although i-\textsc{SCUS} (or equivalently \textsc{SCUS}) is optimal, the returned $F^n\!\left(\bar{P}^n\right)$ is not necessarily concave in $\bar{P}^n$. As a consequence, \textsc{Grad-JSPA} is not guaranteed to converge to a global maximum. 
Nevertheless, we show by numerical results in Section~\ref{sec:num} that it achieves near-optimal WSR performance with low complexity. 

\vspace{-0.2cm}

\subsection{Pseudo-Polynomial Time Optimal Scheme}\label{sec:optJSPA}

The JSPA problem as formulated in~\ref{Psub1} has real variables $\bar{P}^n$ on a continuous search space $\mathcal{F}$. However, the study of NP-hard optimization problems and their approximation requires parameters and variables to be represented by a bounded number of bits~\cite{garey2002computers}, i.e., with bounded precision. This is also a reasonable assumption in practice since MC-NOMA systems are subject to minimum transmit power limitation at the BS and floating-point arithmetic precision of the hardware. As a consequence, we discretize the search space $\mathcal{F}$, in the same way as in~\cite{lei2016power}. Let $\delta$ be the minimum transmit power such that the variables $\bar{P}^n$ can only take value of the form $l\cdot\delta$, for $l\in\{0,1,\ldots,\lfloor\frac{P_{max}}{\delta}\rfloor\}$.
We denote the number of non-zero power values as $J = \lfloor\frac{P_{max}}{\delta}\rfloor$. The feasible set then becomes
\begin{multline*}
\mathcal{F'} \triangleq \{\bar{\bP} \colon\!\sum_{n\in\cN}\bar{P}^n \leq P_{max} \text{ and }0 \leq \bar{P}^n \leq P_{max}^n,\; n\in\cN,
\\ \text{ and } \bar{P}^n = l\cdot\delta,\;l\in\{0,\ldots,J\},\; n\in\cN\}.
\end{multline*}
We rewrite problem~\ref{Psub1} with search space $\mathcal{F'}$ as follows:
\begin{align}
& \underset{\by}{\text{maximize}}
& & \sum_{n\in\cN}\sum_{l=1}^J{c_{n,l}y_{n,l}}, \tag{MCKP}\label{P_MCKP}\\
& \text{subject to}
& & \sum_{n\in\cN}\sum_{l=1}^J{a_{n,l}y_{n,l}}\leq P_{max}, \nonumber\\
&&& \sum_{l=1}^J{a_{n,l}y_{n,l}}\leq P_{max}^n,~n\in\cN, \nonumber\\
&&& \sum_{l=1}^J{y_{n,l}}\leq 1,~n\in\cN, \nonumber\\
&&& y_{n,l} \in \{0,1\},~n\in\cN,~l\in\left[1,J\right], \nonumber
\end{align}
where $c_{n,l} = F^n\!\left(l\cdot\delta\right)$ and $a_{n,l} = l\cdot\delta$. The discretized JSPA problem, denoted by~\ref{P_MCKP}, is known as the multiple choice knapsack problem~\cite{book:knapsack}. It has $N$ disjoint classes each containing $J$ items to be packed into a knapsack of capacity $P_{max}$. Each item has a profit $c_{n,l}$ and a weight $a_{n,l}$, representing respectively the WSR and power consumption of this allocation on subcarrier $n$. The binary variable $y_{n,l}$ takes value $1$ if and only if item $l$ in class $n$ is assigned to the knapsack. The problem consists in assigning at most one item from each class to maximize the sum of their profit. We denote its optimal value by $F^*_{MCKP}$.

As mentioned previously, discretizing~\ref{Psub1} is necessary due to the bounded precision which arises inherently from the study of algorithms and their implementation in practical systems. Besides, \ref{P_MCKP} can be used to approach the continuous solution of~\ref{Psub1} with arbitrary precision. Indeed, Theorem~\ref{th:gap} shows that the discretization error is upper-bounded by a linear function in $\delta$ with a coefficient depending on the system's parameters.
\begin{myth}[Discretization error between $F^*$ and $F^*_{MCKP}$]\label{th:gap}
$ $\newline
The gap between the optimal values of the continuous problem~\ref{Psub1} and its discretized version~\ref{P_MCKP} with step size $\delta$ is upper-bounded by:
\begin{equation*}
F^* - F^*_{MCKP} \leq \delta\sum_{n\in\cN} \max_{k\in\cK}\Bigg\{\frac{W_n w_{\pi_n(k)}}{\left(\bar{P}^{n*}+\tilde{\eta}^n_{\pi_n(k)}\right)\ln\!{(2)}}\Bigg\},
\end{equation*}
where $\bar{P}^{n*}$ is the optimal power budget of~\ref{Psub1} on subcarrier $n\in\cN$. 
\end{myth}
\begin{IEEEproof}
We derive the proof in Appendix~\ref{ap:gap}.
\end{IEEEproof}

The discrete problem~\ref{P_MCKP} can be solved optimally by \textit{dynamic programming by weights} studied in~\cite[Section 11.5]{book:knapsack}.
Based on this idea, we propose \textsc{Opt-JSPA} (see Algorithm~\ref{algoOptJSPA}) to solve~\ref{Psub1}. We first transform~\ref{Psub1} to problem~\ref{P_MCKP}: from line~1 to~5, every item's profit $c_{n,l}$ is computed using i-\textsc{SCUS}. Then, we perform dynamic programming by weights at lines~6-7.
We summarize the optimality and complexity of \textsc{Opt-JSPA} in Theorem~\ref{th:algoOptJSPA}. Detailed analysis of the dynamic programming can be found in Appendix~\ref{ap:algoOptJSPA}.

\begin{algorithm}[ht]
\caption{The pseudo-polynomial time optimal scheme}\label{algoOptJSPA}
\begin{algorithmic}[1]
\MyFunc \textsc{Opt-JSPA}$\left(\left(\mathcal{I}^n\right)_{n\in\cN},M,P_{max},P^n_{max},\delta\right)$
\LineComment Compute the parameters of~\ref{P_MCKP}
\For{$n\in\cN$ and $l\in\left[0,J\right]$}
\State $a_{n,l} \gets l\cdot\delta$
\State $c_{n,l} \gets F^n\left(l\cdot\delta\right)$
\EndFor
\State\Return optimal allocation from the \textit{dynamic program-}
\State \hspace*{10mm} \textit{ming by weights}~\cite{book:knapsack}
\EndFunc
\end{algorithmic}
\end{algorithm}

\begin{table*}[!t]
\vspace{-0.5cm}
\begin{center}
\caption{Comparison of some JSPA schemes proposed in this work and in the literature}
\label{table:JSPA}
\vspace{-0.1cm}
\begin{tabular}{|C{0.23\linewidth}|C{0.28\linewidth}|C{0.41\linewidth}|}
  \hline
  \textbf{Algorithm} & \textbf{Performance guarantee} & \textbf{Complexity for $J$ discrete power values} 
  $\vphantom{\frac{\frac{N}{N}}{\frac{N}{N}}}$\\ \hline
  Monotonic optimization with outer polyblock approximation~\cite{sun2016optimal} & Optimal & Exponential in $K$ and $N$ \\\hline
  \textsc{TSDP}~\cite{lei2016power} & Optimal & $O\!\left(J^2NMK\right) \vphantom{\frac{\frac{N}{N}}{\frac{N}{N}}}$\\\hline
  \textsc{Opt-JSPA} & Optimal  & $O\!\left(NMK^2+JNMK+J^2N\right) 
  \vphantom{\frac{\frac{N}{N}}{\frac{N}{N}}}$\\\hline
  \textsc{\textepsilon-JSPA} & FPTAS, i.e., its performance is within a factor $1-\text{\textepsilon}$ of the optimal, for any $\text{\textepsilon}>0$ &  $O\!\left(NMK^2\!+\!\min\!\Big\{\log\!\left(J\right)\frac{N^2MK}{\text{\textepsilon}}\!+\!\frac{N^3}{\text{\textepsilon}^2}\,,\,JNMK\!+\!J^2N\Big\}\right)$
  $\vphantom{n}$
  \\\hline
  \textsc{Grad-JSPA} & Heuristic & $O\!\left(NMK^2+\log\!\left(J\right)\!NMK\right) \vphantom{\frac{\frac{N}{N}}{\frac{N}{N}}} $\\\hline
\end{tabular}
\end{center}
\vspace{-0.4cm}
\end{table*}

\begin{myth}[Optimality and complexity of \textsc{Opt-JSPA}]\label{th:algoOptJSPA}
$ $\newline
Given a minimum transmit power $\delta$, algorithm \textsc{Opt-JSPA} computes the optimal of~\ref{Psub1} on the discrete set $\mathcal{F'}$. Its computational complexity is $O\!\left(NMK^2+JNMK+J^2N\right)$, which is pseudo-polynomial in $J$.
\end{myth}
\begin{IEEEproof}
We explain the principle of dynamic programming by weights and derive its complexity in Appendix~\ref{ap:algoOptJSPA}.
\end{IEEEproof}
\textsc{Opt-JSPA} is said to be pseudo-polynomial since it depends on the total number of power values $J$, whereas all system's parameters and variables are encoded in $O\!\left(\log\!\left(J\right)\right)$ bits. As a consequence, in practical systems, the contribution of $J$ to the computation time is way higher than parameters $N$, $K$, $M$. 

\subsection{Fully Polynomial-Time Approximation Scheme}\label{sec:epsilonJSPA}

We develop a FPTAS to avoid the pseudo-polynomial complexity in $J$ that is inherent to the optimal schemes \textsc{Opt-JSPA} and \textsc{TSDP}~\cite{lei2016power}. 
According to~\cite{vazirani2013approximation}, an algorithm is said to be a FPTAS if it outputs a solution within a factor $1-\text{\textepsilon}$ of the optimal, for any $\text{\textepsilon}>0$. Moreover, its running time is bounded by a polynomial in both the input size and $1/\text{\textepsilon}$. 
A FPTAS is the best trade-off one can hope for an NP-hard optimization problem in terms of performance guarantee and complexity, assuming P $\neq$ NP. 

The proposed FPTAS, called \textsc{\textepsilon-JSPA} (see Algorithm~\ref{algoepsilonJSPA}), is based on dynamic programming with scaled profits. Scaling the profits is a common technique to reduce the number of items computed in~\ref{P_MCKP}. First, we compute an estimation $U$ of~\ref{P_MCKP}'s optimal value, such that $U \geq F^*_{MCKP} \geq U/4$. We explain the estimation procedure in Appendix~\ref{ap:estimationU}. Then, instead of computing all $JN$ profit values $c_{n,l}$, we only consider the subset $L_n$ of items on each subcarrier $n$ such that:
\begin{equation*}
L_n \triangleq \{l'\leq J, l \leq \frac{4N}{\text{\textepsilon}}-1\colon\! c_{n,l'} \geq l\frac{\text{\textepsilon}U}{4N} > c_{n,l'-1}\}.
\end{equation*}
This can be seen as considering only one profit value per interval of the form $\left[(l-1)\cdot{\text{\textepsilon}U}/{4N},l\cdot{\text{\textepsilon}U}/{4N}\right]$, for $l\in\{1,\ldots,4N/\text{\textepsilon}\}$.
Each $L_n$, for $n\in\cN$, can be obtained by multi-key binary search~\cite{tarek2008multi}. All function evaluations required by the multi-key binary search are done by i-\textsc{SCUS}.

Finally, we apply the \textit{dynamic programming by profits}~\cite[Section 11.5]{book:knapsack} in lines~5-6. It is known that the optimal solution obtained by dynamic programming by profits considering only items in $L_n$, differs from $F^*_{MCKP}$ by at most a factor $1-\text{\textepsilon}$. The performance of \textsc{\textepsilon-JSPA} are summarized in Theorem~\ref{th:algoEpsJSPA}. We provide more details on the estimation $U$ in Appendix~\ref{ap:estimationU} and the dynamic programming by profits in Appendix~\ref{ap:algoEpsJSPA}.

\begin{myth}[Performance and complexity of \textsc{\textepsilon-JSPA}]\label{th:algoEpsJSPA}
$ $\newline
Given a minimum transmit power $\delta$ and an approximation factor \textepsilon, algorithm \textsc{\textepsilon-JSPA} computes an \textepsilon-approximation of~\ref{Psub1} on the discrete set $\mathcal{F'}$. The algorithm is a FPTAS with asymptotic complexity: $O\!\left(NMK^2\!+\!\min\!\Big\{\log\!\left(J\right)\frac{N^2MK}{\text{\textepsilon}}\!+\!\frac{N^3}{\text{\textepsilon}^2}\,,\,JNMK\!+\!J^2N\Big\}\right)$.
\end{myth}
\begin{IEEEproof}
We derive this result in Appendix~\ref{ap:algoEpsJSPA}, using the dynamic programming by profits and the estimation procedure studied in Appendix~\ref{ap:estimationU}.
\end{IEEEproof}

\begin{algorithm}[ht]
\caption{The proposed FPTAS (\textsc{\textepsilon-JSPA})}\label{algoepsilonJSPA}
\begin{algorithmic}[1]
\MyFunc \textsc{\textepsilon-JSPA}$\left(\left(\mathcal{I}^n\right)_{n\in\cN},M,P_{max},P^n_{max},\delta,\text{\textepsilon}\right)$
\State Compute an estimation $U$ of $F^*_{MCKP}$
\For{$n\in\cN$}
\State Get $a_{n,l}$, $c_{n,l}$, for $l\in L_n$ by multi-key binary search
\EndFor
\State\Return \textepsilon-approximate allocation from the \textit{dynamic}
\State \hspace*{10mm} \textit{programming by profits}~\cite{book:knapsack}
\EndFunc
\end{algorithmic}
\end{algorithm}

\vspace{-0.6cm}

\subsection{Comparison of JSPA Algorithms}
In Table~\ref{table:JSPA}, we compare the performance and complexity of the proposed algorithms with JSPA schemes in the literature. Reference~\cite{sun2016optimal} studied an optimal monotonic optimization framework, which has exponential complexity in $K$ and $N$. The two-stage dynamic programming algorithm (\textsc{TSDP}) proposed by Lei et al. has complexity $O\!\left(J^2NMK\right)$ according to~\cite[Theorem 13]{lei2016power}. Both \textsc{TSDP} and the proposed \textsc{Opt-JSPA} are optimal. However, \textsc{Opt-JSPA} has lower complexity than \textsc{TSDP}. Indeed, the right term $J^2N$ is lower by a factor $MK$, the middle term $JNMK$ by a factor $J$. The left term $NMK^2$ also improves the complexity, since reference~\cite{fu2018subcarrier} shows that in practical systems $J = \Theta\left(\min\{K,MN\}\right)$. This result is verified by simulation in Section~\ref{sec:num}.
\textsc{\textepsilon-JSPA} is the proposed FTPAS. Its complexity is bounded by a polynomial in $N/\text{\textepsilon}$ and $\log\!\left(J\right)$. If $N/\text{\textepsilon} = O(J)$, it has lower complexity than \textsc{Opt-JSPA}. Otherwise, if $N/\text{\textepsilon} = \Omega(J)$, its complexity is asymptotically equivalent to \textsc{Opt-JSPA}'s complexity. This means that for very low error rate \textepsilon, the complexity of \textsc{\textepsilon-JSPA} tends to that of \textsc{Opt-JSPA}.
Finally, \textsc{Grad-JSPA} is a heuristic. Its performance is evaluated through simulation in the next section. When applied in a discrete setting, the error tolerance or precision $\xi$ is related to $\delta = 2\xi$. Hence, its complexity is proportional to $\log\!\left(J\right)$, which is way lower than the optimal schemes with pseudo-polynomial complexity due to $J$.

\vspace{-0.2cm}

\section{Numerical Results}\label{sec:num}
\begin{figure}[!t]
\vspace{-0.7cm}
\centering
\includegraphics[width=0.93\linewidth]{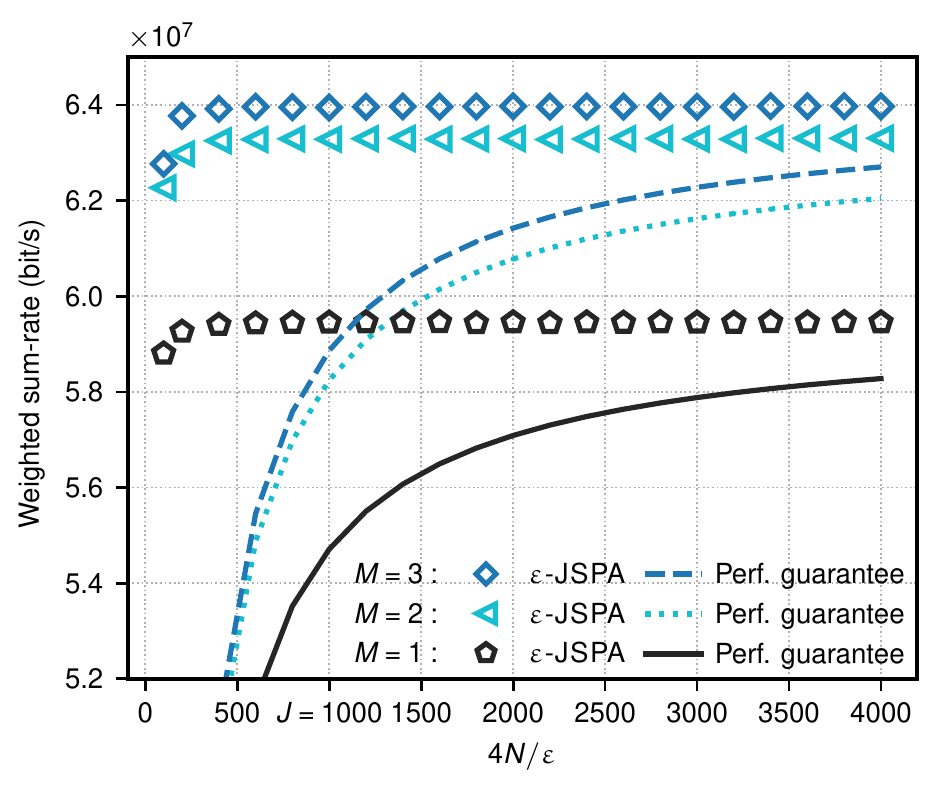}
\vspace{-0.3cm}
\caption{WSR of \textsc{\textepsilon-JSPA} and its guaranteed performance bound versus $4N/\text{\textepsilon}$}
\label{fig_wsr_epsilon}
\end{figure}

\begin{figure}[!t]
\vspace{-0.4cm}
\centering
\includegraphics[width=0.95\linewidth]{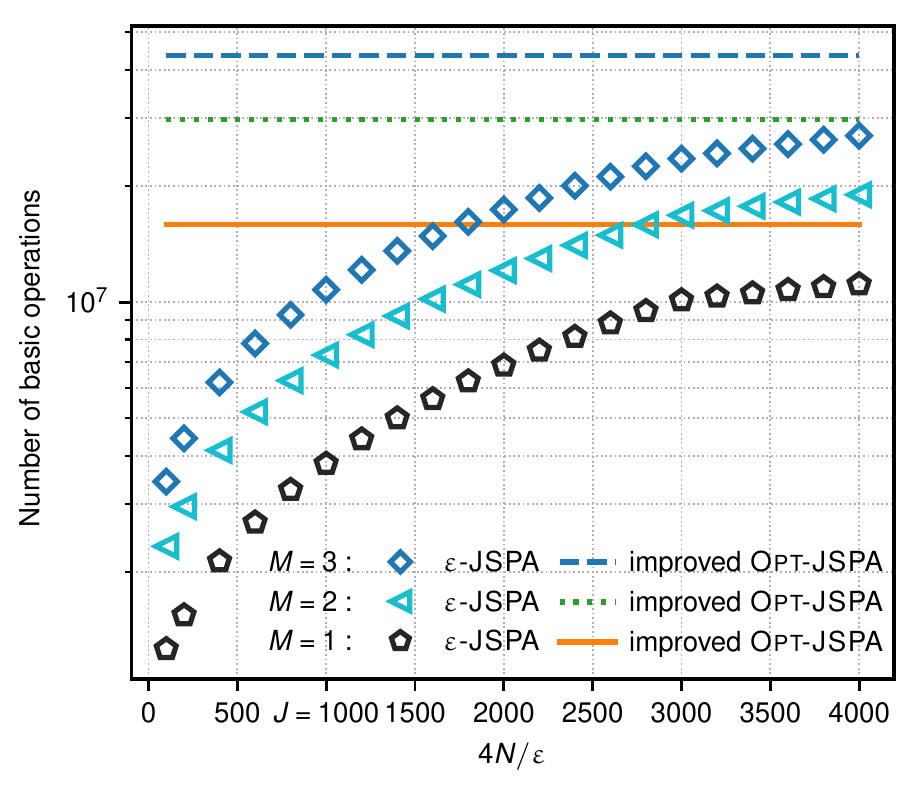}
\vspace{-0.3cm}
\caption{Number of basic operations performed by \textsc{\textepsilon-JSPA} versus $4N/\text{\textepsilon}$}
\label{fig_cplx_epsilon}
\end{figure}

We evaluate the WSR and computational complexity of \textsc{Opt-JSPA}, \textsc{\textepsilon-JSPA} and \textsc{Grad-JSPA} through numerical simulations. We compare them with the optimal benchmark scheme \textsc{TSDP} introduced in~\cite{lei2016power}.
We consider a hexagonal cell of diameter $1000$ meters, with one BS located at its center and $K$ users distributed uniformly at random in the cell. The users' weights are generated uniformly at random in $\left[0,1\right]$. The number of users $K$ varies from $5$ to $60$, and the number of subcarriers is $N = 20$. We assume a system bandwidth of $W=5$~MHz and $W_n=W/N$ for all $n$. We follow the radio propagation model of~\cite{Greentouch}, including path loss, shadowing and Rayleigh fading. The minimum transmit power is $\delta = 0.01\,\text{W}$. The cellular power budget is $P_{max} = 10\,\text{W}$, therefore the number of power values is $J=10^3$.
Each point in the following figures is the average value obtained over $1000$ random instances. Only Fig.~\ref{fig_wsr_epsilon} and~\ref{fig_cplx_epsilon} represent a single instance. The simulation parameters and channel model are summarized in Table~\ref{parameter}. 
\begin{table}[!t]
\vspace{-0.2cm}
\begin{center}
\caption{Simulation parameters}
\label{parameter}
\scalebox{0.95}{
\begin{tabular}{|c|c|}
  \hline
  \textbf{Parameter}  & \textbf{Value} \\ \hline
  Cell radius  & $1000$ m \\ \hline
  Min. distance from user to BS & $35$ m \\ \hline
  Carrier frequency & $2$ GHz \\ \hline
  Path loss model & $128.1 + 37.6 \log_{10} d$ dB, $d$ is in km \\ \hline
  Shadowing  & Log-normal, $10$ dB standard deviation \\ \hline
  Fading & Rayleigh fading with variance 1 \\ \hline
  Noise power spectral density & $-174$ dBm/Hz \\ \hline
  System bandwidth $W$ & $5$ MHz \\ \hline
  Number of subcarriers $N$ & $20$  \\ \hline
  Number of users $K$ & $5$ to $60$  \\ \hline
  Total power budget $P_{max}$ & $10$ W \\ \hline
  Minimum transmit power $\delta$ & $0.01$ W\\ \hline
  Number of power values $J$ & $10^3$ \\ \hline
  Error tolerance $\xi$ & $10^{-4}$\\ \hline
  Parameter $M$ & $1$ (OMA), $2$ and $3$ (NOMA)\\
  \hline
\end{tabular}}
\end{center}
\vspace{-0.2cm}
\end{table}

\begin{figure}[!t]
\vspace{-0.3cm}
\centering
\includegraphics[trim={0 0 0 0.2cm},clip,width=0.95\linewidth]{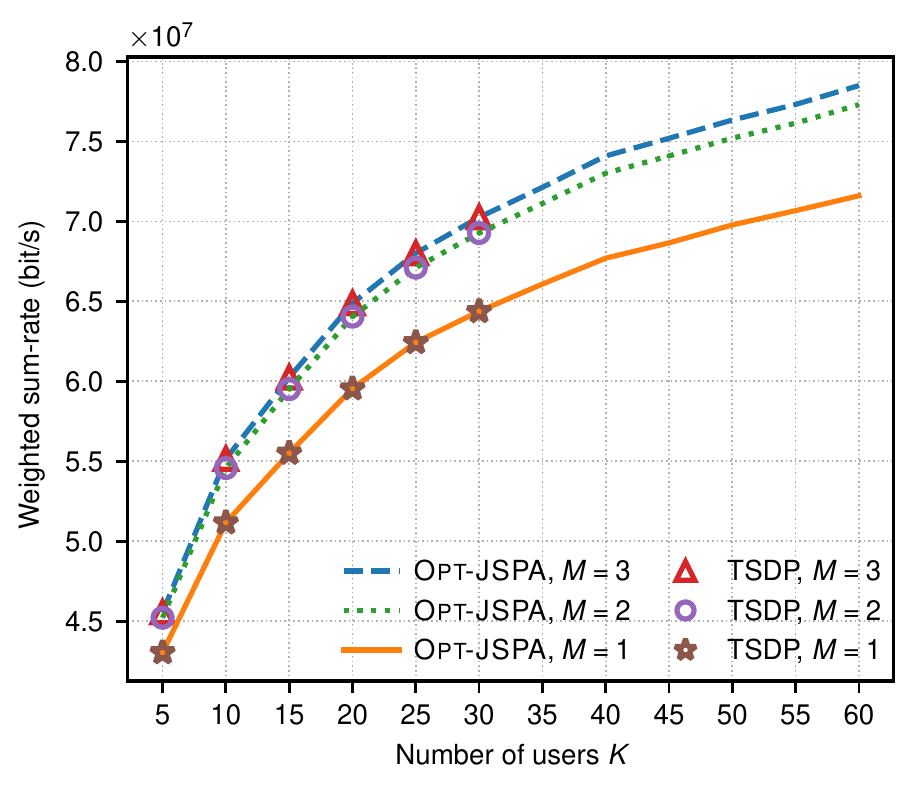}
\vspace{-0.4cm}
\caption{WSR of the optimal schemes for different number of users $K$}
\label{fig_wsr}
\end{figure}

\begin{figure}[!t]
\vspace{-0.4cm}
\centering
\includegraphics[width=0.92\linewidth]{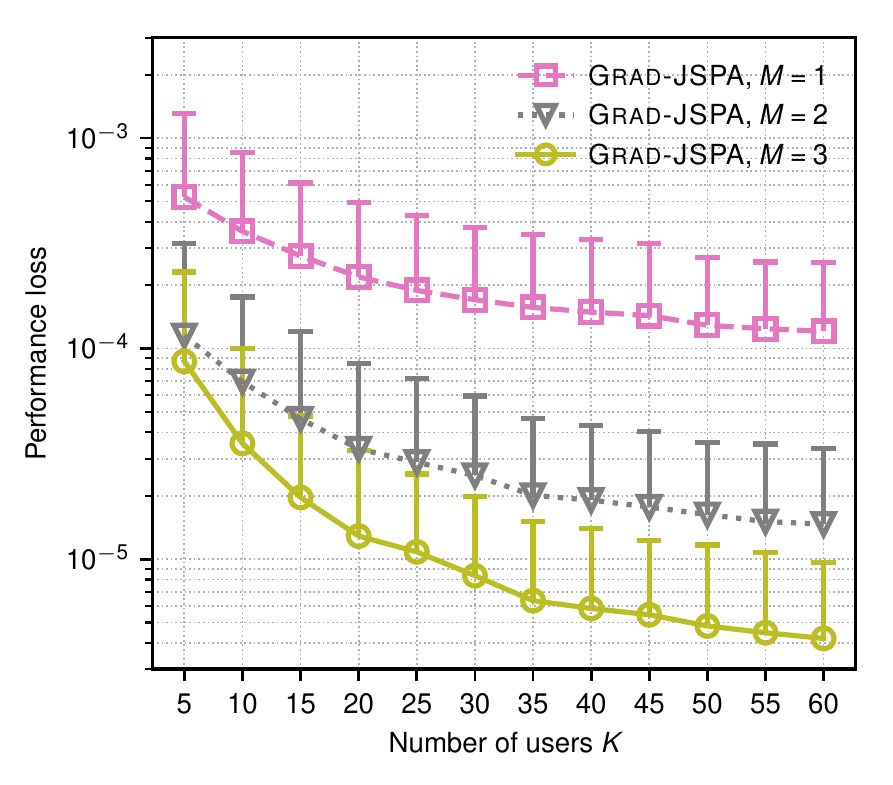}
\vspace{-0.4cm}
\caption{Performance loss of \textsc{Grad-JSPA} compared to the optimal WSR}
\label{fig_perf_loss}
\end{figure}

\begin{figure}[!t]
\vspace{-0.3cm}
\centering
\includegraphics[trim={0.1cm 0 0.1cm 0.1cm},clip,width=0.96\linewidth]{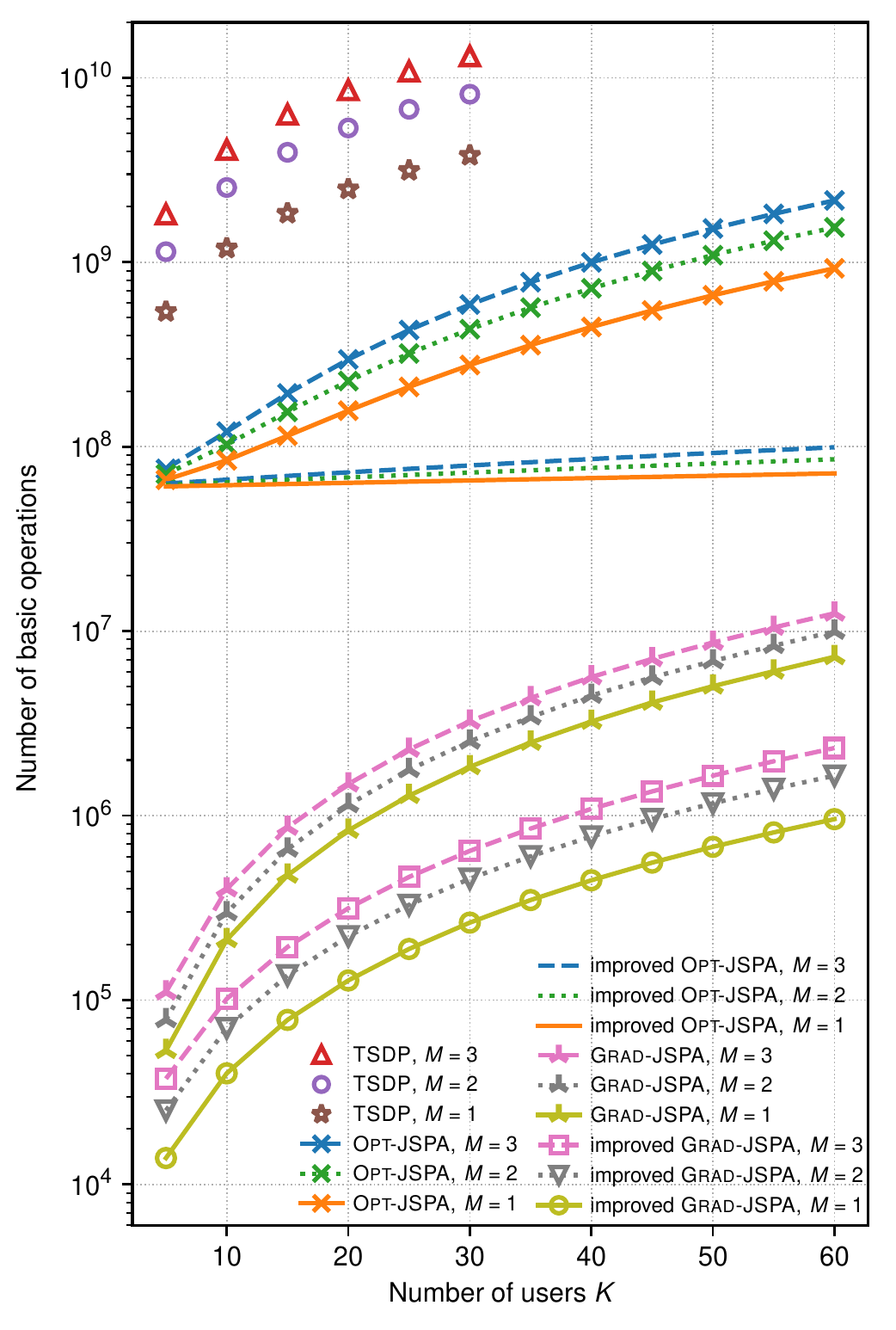}
\vspace{-0.4cm}
\caption{Number of basic operations performed by each algorithm versus $K$}
\label{fig_cplx}
\end{figure}

Fig.~\ref{fig_wsr} shows the WSR performance of \textsc{Opt-JSPA} and \textsc{TSDP}, for $M = 1$, $2$ and $3$. We only simulate \textsc{TSDP} for $K = 5$ to $30$, due to its high running time complexity. We see that \textsc{Opt-JSPA} and \textsc{TSDP} achieve the same WSR performance, which is consistent with the fact that they are both optimal. Indeed, the optimality of \textsc{Opt-JSPA} is shown in Theorem~\ref{th:algoOptJSPA}, and the optimality of \textsc{TSDP} has been proven in~\cite[Theorem 13]{lei2016power}. Although both algorithms have the same performance, we will see further on in Fig.~\ref{fig_cplx} that \textsc{Opt-JSPA} has lower computational complexity than \textsc{TSDP}. The performance gain of NOMA with $M = 2$ (resp. $M=3$) over OMA (i.e., $M=1$) is about $8\%$ (resp. $10\%$), for $K=60$.

Fig.~\ref{fig_perf_loss} illustrates the performance loss of \textsc{Grad-JSPA} compared to the optimal, for $M = 1$, $2$ and $3$. The performance loss is defined as: 
\begin{equation*}
\frac{\text{Optimal WSR} - \text{\textsc{Grad-JSPA} WSR}}{\text{Optimal WSR}}.
\end{equation*} 
The markers represent the average performance loss, while the upper intervals indicate the 90th percentile. For example, for $K=10$ and $M=1$, $90\%$ of \textsc{Grad-JSPA} results have less than $9\times 10^{-4}$ of performance loss. We observe that the average performance loss is always below $6\times 10^{-4}$. Hence, our proposed heuristic \textsc{Grad-JSPA} achieves near-optimal solutions in these simulation settings.
It is also suitable for large systems, since the performance loss decreases as $K$ or $M$ increases. 

In Fig.~\ref{fig_cplx}, we count the number of basic operations
(additions, multiplications, comparisons) performed by each
algorithm, which reflects their computational complexity. The term ``improved" in the legend represents the complexity of \textsc{Opt-JSPA} and \textsc{Grad-JSPA} when using i-\textsc{SCPC} and i-\textsc{SCUS} instead of \textsc{SCPC} and \textsc{SCUS}. There is a significant speed up by employing i-\textsc{SCPC} and i-\textsc{SCUS} as basic building blocks. Indeed, for $K=60$ and $M=1$, $2$ or $3$, there is a factor of at least $10$ between \textsc{Opt-JSPA} and its improved version. Besides, the improved \textsc{Opt-JSPA} outperforms \textsc{TSDP} in terms of complexity. For instance, \textsc{Opt-JSPA} reduces the complexity by a factor $330$, for $K=30$ and $M=3$. Finally, \textsc{Grad-JSPA} has low complexity, which makes it a good choice for practical implementation. 

Fig.~\ref{fig_wsr_epsilon} and~\ref{fig_cplx_epsilon} present the WSR and complexity of \textsc{\textepsilon-JSPA} versus $4N/\text{\textepsilon}$. We choose such a normalized x-axis, as it is equal to the number of items evaluated in each subcarrier, i.e., $|L_n| = 4N/\text{\textepsilon}$. It can be directly compared to $J$, which is the total number of items in each subcarrier in the discretized problem~\ref{P_MCKP}. Here, we simulate a single instance with $K=60$ users to show how \textsc{\textepsilon-JSPA} behaves as a function of $\text{\textepsilon}$. In Fig.~\ref{fig_wsr_epsilon}, we also present its performance guarantee. Recall that the performance guarantee is $1-\text{\textepsilon}$ times the optimal. As expected, \textsc{\textepsilon-JSPA} is always above its performance guarantee. As $N/\text{\textepsilon}$ increases, the approximation guarantee tends to the optimal. In this instance, the algorithm already achieves a near-optimal solution for $4N/\text{\textepsilon}=400$, i.e., $\text{\textepsilon}=0.2$.
In Fig.~\ref{fig_cplx_epsilon}, we also plot the complexity of the improved \textsc{Opt-JSPA} for comparison. As explained in Section~\ref{sec:epsilonJSPA}, the complexity increases with $N/\text{\textepsilon}$ and becomes (asymptotically) equal to that of \textsc{Opt-JSPA} for $N/\text{\textepsilon} = \Omega(J)$. In this regime, there is apparently no benefit of using \textsc{\textepsilon-JSPA}, since \textsc{Opt-JSPA} achieves the optimal with the same complexity. 
Nevertheless, in practice, we can see that even for $4N/\text{\textepsilon}\geq J$, \textsc{\textepsilon-JSPA} has less operations than \textsc{Opt-JSPA}. This is because the number of items computed by \textsc{\textepsilon-JSPA} increases slowly and smoothly as a function of \textepsilon. This behavior is not captured in the asymptotic complexity (big-$O$ notation). This is verified in Fig.~\ref{fig_cplx_epsilon} for up to $4J = 4000$.
In summary, \textsc{\textepsilon-JSPA} allows us to control the trade-off between WSR and complexity with $\text{\textepsilon}$.

\section{Discussion on Possible Generalizations}\label{sec:discussion}
In this work, we assume that the channel gains are perfectly known. Two more realistic models using only partial CSI can be considered instead: \textit{imperfect CSI} studied in~\cite{wei2017optimal,choi2017joint}, for which the channel gains are given with a known estimation error probability distribution, and \textit{second order statistics} (SOS) adopted in~\cite{yang2015performance}, for which only the distances between users and BS are known. 
We believe that our framework can be extended to these cases by maximizing the expected WSR depending on stochastic channel gains, while the power constraints remain unchanged (i.e., non-stochastic). The challenge would be to characterize such a stochastic objective function in many scenarios. This is a possible future research direction.

As multi-antenna technologies are becoming more and more important in 5G and Beyond 5G systems~\cite{dahlman20185g,wong2017key}, it would be interesting to extend the current work to multi-antenna transmissions. Paper~\cite{islam2018resource} states that MC-NOMA with multiple antennas is a much more complex problem which requires to develop novel low complexity solutions. One may draw inspiration from the work of Sun et al.~\cite{sun2017optimalMISO} which generalizes the monotonic optimization framework of~\cite{sun2016optimal} to a MC-MISO-NOMA system. A similar approach might be adopted in our framework, which is explained as follows: Since the SIC decoding order in multi-antenna MC-NOMA systems does not only depend on the channel gains but also on the beamforming (BF), the user clustering (UC) and BF have to be jointly optimized to achieve optimal or approximate performance. Hence, the idea would be to extend \textsc{SCUS} to a joint UC and BF optimization scheme. This scheme can then be integrated in \textsc{Opt-JSPA}, \textsc{\textepsilon-JSPA} and \textsc{Grad-JSPA}, while preserving their performance guarantees. Although optimal joint UC and BF remains a difficult open problem, existing heuristics can be adopted instead. For example, the schemes of paper~\cite{ali2016non} have shown to outperform classical (OMA-based) MIMO systems and other multi-antenna NOMA algorithms by simulations.


\section{Conclusion}\label{sec:ccl}
In this work, we investigate the WSR maximization in MC-NOMA with cellular power constraint. We improve the complexity of the single-carrier power control (\textsc{SCPC}) and user selection (\textsc{SCUS}) procedures using precomputation. These improved schemes are denoted by i-\textsc{SCPC} and i-\textsc{SCUS}. We develop three algorithms to solve the JSPA problem, based on i-\textsc{SCPC} and i-\textsc{SCUS}. Firstly, \textsc{Opt-JSPA} gives optimal results with lower complexity than current state-of-the-art optimal schemes, i.e., TSDP~\cite{lei2016power} and monotonic optimization~\cite{sun2016optimal}. Secondly, \textsc{\textepsilon-JSPA} is a FPTAS. It achieves a controllable and tight trade-off between approximation guarantee and complexity. \textsc{Opt-JSPA} and \textsc{\textepsilon-JSPA} are both suitable for performance benchmarking. Finally, \textsc{Grad-JSPA} is a heuristic. We show by simulation that it has near-optimal WSR with low and practical complexity. 

\bibliographystyle{IEEEtran}
\bibliography{IEEEabrv,reference}

\vspace{-0.2cm}
\appendix

We first provide in Lemma~\ref{le:propOptimal_sumf} an important property on the solution maximizing $\sum_{l=1}^i{\tilde{f}_{l,l}^n}$ subject to $C2'\text{--}3'$, for $i\leq |\cU'_n|$. This Lemma will be used in Appendices~\ref{ap:thalgoSCPC} and~\ref{ap:derivative}.

\begin{myle}\label{le:propOptimal_sumf}
$ $\newline
Assume we are given a subcarrier $n\in\cN$, a set $\cU'_n$ of active users, a power budget $\bar{P}^n$, and an index $i$. Let $x_{1_n}^n,\ldots,x_{i_n}^n$ be the allocation maximizing $\sum_{l=1}^i{\tilde{f}_{l,l}^n\left(\cU'_n,x_{l_n}^n\right)}$, while also satisfying $C2'\text{--}3'$, i.e., $\bar{P}^n \geq x_{1_n}^n \geq\cdots\geq x_{i_n}^n \geq 0$. 
$x_{1_n}^n,\ldots,x_{i_n}^n$ can be partitioned into sequences of consecutive terms with the same value. That is, sequences of the form $x_{q_n}^n,\ldots,x_{q'_n}^n$, where $x_{q_n}^n=\cdots=x_{q'_n}^n$ and $1\leq q\leq q'\leq |\cU'_n|$, $q = 1$ or $x_{(q-1)_n}^n > x_{q_n}^n$, $q' = |\cU'_n|$ or $x_{q'_n}^n < x_{(q'+1)_n}^n$.
Any such sequence satisfies:
\begin{equation*}
x_{q_n}^n=\cdots=x_{q'_n}^n=\textsc{Argmax}\tilde{f}\!\left(q,q',\mathcal{I}^n,\cU'_n,\bar{P}^n\right).
\end{equation*}
\end{myle}
\begin{IEEEproof}
In this proof, we simplify notation $\tilde{f}_{l,l}^n\left(\cU'_n,\cdot\right)$ as $\tilde{f}_{l,l}^n\left(\cdot\right)$.
Let $x_{q_n}^n,\ldots,x_{q'_n}^n$ be a sequence of consecutive terms with the same value, as defined in Lemma~\ref{le:propOptimal_sumf}. Assume, for the sake of contradiction, that $x_{q_n}^n=\cdots=x_{q'_n}^n\neq x^*$, where $x^* = \textsc{Argmax}\tilde{f}\!\left(q,q',\mathcal{I}^n,\cU'_n,\bar{P}^n\right)$. Without loss of generality, we consider the case $x_{q_n}^n < \textsc{Argmax}\tilde{f}\!\left(q,q',\mathcal{I}^n,\cU'_n,\bar{P}^n\right)$ and $q > 1$. Let $y_{1_n}^n,\ldots,y_{i_n}^n$ be an allocation defined as:
\begin{equation}\label{ap:le:proof11_cases}
y_{l_n}^n \triangleq \begin{cases}
    \min\{x_{(q-1)_n}^n,x^*\}, & \text{if $q \leq l \leq q'$},\\
    x_{l_n}^n, & \text{otherwise}.
  \end{cases} 
\end{equation}
We have the following inequalities:
\begin{align}
&\sum_{l=1}^i{\tilde{f}_{l,l}^n\left(y_{l_n}^n\right)} = \sum_{l\notin\{q,\ldots,q'\}}{\tilde{f}_{l,l}^n\left(x_{l_n}^n\right)}+\tilde{f}_{q,q'}^n\left(y_{l_n}^n\right), \label{ap:le:proof11_ineq1}\\
&\quad> \sum_{l\notin\{q,\ldots,q'\}}{\tilde{f}_{l,l}^n\left(x_{l_n}^n\right)}+\tilde{f}_{q,q'}^n\left(x_{l_n}^n\right)= \sum_{l=1}^i{\tilde{f}_{l,l}^n\left(x_{l_n}^n\right)}. \label{ap:le:proof11_ineq2}
\end{align}
Equality~\eqref{ap:le:proof11_ineq1} comes from the definition in~\eqref{ap:le:proof11_cases}. According to Lemma~\ref{le:maxfij}, $\tilde{f}_{q,q'}^n$ is increasing on $[0,x^*]$, which implies inequality~\eqref{ap:le:proof11_ineq2}. In summary, $y_{1_n}^n,\ldots,y_{i_n}^n$ satisfies $C2'\text{--}3'$ by its definition in~\eqref{ap:le:proof11_cases}, and it achieves greater value of $\sum_{l=1}^i{\tilde{f}_{l,l}^n}$ than $x_{1_n}^n,\ldots,x_{i_n}^n$. This is a contradiction, therefore it must be that \begin{equation}\label{ap:le:proof11_finalEq1}
x_{q_n}^n \geq \textsc{Argmax}\tilde{f}\!\left(q,q',\mathcal{I}^n,\cU'_n,\bar{P}^n\right).
\end{equation}
If $q=1$, the same reasoning can be applied by replacing $\min\{x_{(q-1)_n}^n,x^*\}$ by $\bar{P}^n$ in Eqn.~\eqref{ap:le:proof11_cases}. 
We can perform a similar proof by contradiction on the case $x_{q_n}^n > \textsc{Argmax}\tilde{f}\!\left(q,q',\mathcal{I}^n,\cU'_n,\bar{P}^n\right)$ to deduce that:
\begin{equation}\label{ap:le:proof11_finalEq2}
x_{q_n}^n \leq \textsc{Argmax}\tilde{f}\!\left(q,q',\mathcal{I}^n,\cU'_n,\bar{P}^n\right).
\end{equation}
The desired result follows from~\eqref{ap:le:proof11_finalEq1} and~\eqref{ap:le:proof11_finalEq2}.
\end{IEEEproof}

\subsection{Proof of Lemma~\ref{le:equivalence}}\label{ap:transf}
The objective of~\ref{P} can be written as:
\begin{align*}
&\sum_{k\in\cK}{w_k\sum_{n\in\cN}{R_k^n\left(\bp^n\right)}} = \sum_{n\in\cN}\sum_{k\in\cK}{w_k R_k^n\left(\bp^n\right)},\label{eq:rewrite1}\\
&\stackrel{\text{(b)}}{=} \sum_{n\in\cN}W_n\sum_{i=1}^{K}{w_{\pi_n(i)}\log_2\left(\frac{\sum_{j=i}^{K}{p_{\pi_n(j)}^n}+\tilde{\eta}^n_{\pi_n(i)}}{\sum_{j=i+1}^{K}{p_{\pi_n(j)}^n}+\tilde{\eta}^n_{\pi_n(i)}}\right)},\\
&\stackrel{\text{(c)}}{=} \sum_{n\in\cN}W_n\sum_{i=1}^{K}{\log_2\left(\frac{\left(\sum_{j=i}^{K}{p_{\pi_n(j)}^n}+\tilde{\eta}^n_{\pi_n(i)}\right)^{w_{\pi_n(i)}}}{\left(\sum_{j=i+1}^{K}{p_{\pi_n(j)}^n}+\tilde{\eta}^n_{\pi_n(i)}\right)^{w_{\pi_n(i)}}}\right)},\\
&\stackrel{\text{(d)}}{=} \sum_{n\in\cN}W_n{\left[w_{\pi_n(1)}\log_2\left(\sum_{j=1}^{K}{p_{\pi_n(j)}^n}+\tilde{\eta}^n_{\pi_n(1)}\right)\right.}\\
&\phantom{=} + \sum_{i=2}^{K}{\log_2\left(\frac{\left(\sum_{j=i}^{K}{p_{\pi_n(j)}^n}+\tilde{\eta}^n_{\pi_n(i)}\right)^{w_{\pi_n(i)}}}{\left(\sum_{j=i}^{K}{p_{\pi_n(j)}^n}+\tilde{\eta}^n_{\pi_n(i-1)}\right)^{w_{\pi_n(i-1)}}}\right)} \\
&\phantom{=} + \left. w_{\pi_n(K)}\log_2\left(\frac{1}{\tilde{\eta}^n_{\pi_n(K)}}\right) \right].
\end{align*}
Equality~(b) comes from the definition in~\eqref{eq:R_2}. At~(c), the weights $w_{\pi_n(i)}$ are put inside the logarithm. Finally, (d) is obtained by combining the numerator of the $i$-th term with the denominator of the $(i-1)$-th term, for $i\in\{2,\ldots,K\}$. 

By applying the change of variables shown in~\eqref{eq:change}, we derive the equivalent problem~\ref{P2}. The constant term is $A = \sum_{n\in\cN}w_{\pi_n(K)}\log_2\left(1/\tilde{\eta}^n_{\pi_n(K)}\right)$.
Constraints $C1'$ and $C2'$ are respectively equivalent to $C1$ and $C2$ since $x_1^n = \sum_{j=1}^{K}{p_{\pi_n(j)}^n} = \sum_{k\in\cK}{p_k^n}$, for $n\in\cN$. Constraints $C3'$ and $C3''$ come from $C3$ and the fact that $x_i^n - x_{i+1}^n = p_{\pi_n(i)}^n$, for any $i\in\{1,\ldots,K\}$ and $n\in\cN$. In the same way, the active users set in $C4'$ is defined as ${\cU'_n \triangleq \{i\in\{1,\ldots,K\} \colon x_i^n > x_{i+1}^n\}}$. 
\hfill\IEEEQEDhere

\subsection{Proof of Lemma~\ref{le:maxfij}}\label{ap:propf}
We study the first and second derivatives of $f_{j,i}^n$, denoted by ${f_{j,i}^n}'$ and ${f_{j,i}^n}''$. If $j=1$, then we have:
\begin{equation}\label{proof:first_deriv_f}
{f_{1,i}^n}'(x) = \frac{W_n w_{\pi_n(i)}}{\left(x+\tilde{\eta}^n_{\pi_n(i)}\right)\ln\!{(2)}},
\end{equation}
which is strictly positive and decreasing for $x\geq 0$. Hence, $f_{1,i}^n$ is increasing and concave. 
For $j > 1$, the first and second derivatives are as follows:
\begin{align*}
{f_{j,i}^n}'(x) &= \frac{W_n}{\ln\!{(2)}}\!\left(\frac{ w_{\pi_n(i)}}{x+\tilde{\eta}^n_{\pi_n(i)}}-\frac{ w_{\pi_n(j-1)}}{x+\tilde{\eta}^n_{\pi_n(j-1)}}\right),\\
{f_{j,i}^n}''(x) &= \frac{W_n}{\ln\!{(2)}}\!\left(\frac{ w_{\pi_n(j-1)}}{\left(x+\tilde{\eta}^n_{\pi_n(j-1)}\right)^2}-\frac{ w_{\pi_n(i)}}{\left(x+\tilde{\eta}^n_{\pi_n(i)}\right)^2}\right).
\end{align*}
We know that $\tilde{\eta}^n_{\pi_n(j-1)} \geq \tilde{\eta}^n_{\pi_n(i)}$ by construction of the optimal decoding order in Eqn.~\eqref{SIC_DL}. If, in addition, we have $w_{\pi_n(i)} \geq w_{\pi_n(j-1)}$, then ${f_{j,i}^n}'(x)\geq 0$ and ${f_{j,i}^n}''(x)\leq 0$ for all $x\geq 0$. We deduce that 
$f_{j,i}^n$ is increasing and concave. This proves the first point of Lemma~\ref{le:maxfij}. Now suppose that $w_{\pi_n(i)} < w_{\pi_n(j-1)}$ instead. Values $c_1$ and $c_2$ defined in Lemma~\ref{le:maxfij} are the unique roots of the first and second derivatives, i.e., ${f_{j,i}^n}'(c_1) = 0$ and ${f_{j,i}^n}''(c_2) = 0$. ${f_{j,i}^n}'$ is positive on $\left(-\tilde{\eta}_{\pi_n(j-1)},c_1\right)$ and negative on $\left(c_1,\infty\right)$. This implies that $f_{j,i}^n$ is unimodal and has a unique global maximum at $c_1$ for $x>0$. Similarly, ${f_{j,i}^n}''$ is negative on $\left(-\tilde{\eta}_{\pi_n(j-1)},c_2\right)$ and positive on $\left(c_2,\infty\right)$. Therefore, $f_{j,i}^n$ is concave before $c_2$ and convex after $c_2$. This proves the second point of Lemma~\ref{le:maxfij}.
\hfill\IEEEQEDhere

\subsection{Proof of Theorem~\ref{th:algoSCPC}}\label{ap:thalgoSCPC}

The complexity and optimality proofs of \textsc{SCPC} are presented below.

\textbf{\textit{Complexity analysis:}} At each for loop iteration $i$, the while loop at line 6 has at most $i$ iterations. Thus, the worst case complexity is proportional to $\sum_{i=1}^{|\cU'_n|}{i} = O\!\left(|\cU'_n|^2\right)=O\!\left(M^2\right)$.

\textbf{\textit{Optimality analysis:}} Without loss of generality, we can suppose that the $x_{i_n}^n$'s are initialized to zero. We will prove by induction that at the end of each iteration $i$ at line~10 of Algorithm~\ref{algoSCPC}, the following loop invariants are true:
\begin{enumerate}[label={$H_\arabic*(i)$:}, ref={$H_\arabic*$}, leftmargin=*, itemindent=3em]
\item ${\sum_{l=1}^{i}{\tilde{f}_{l,l}^n}}$ is maximized by $x_{1_n}^n,\ldots,x_{i_n}^n$,\label{h1}
\item $C2'\text{--}3'$ is satisfied, i.e., $\bar{P}^n \geq x_{1_n}^n \geq \cdots \geq x_{i_n}^n \geq 0$.\label{h2}
\end{enumerate}
\textbf{\textit{Basis:}} For $i=1$, $x^*$ computed at line~3 is indeed the optimal of $\tilde{f}_{1,1}^n$. The while loop has no effect since $j = 0 < 1$, therefore $x_{1_n}^n \gets x^*$ and statements~\ref{h1}$(1)$ and~\ref{h2}$(1)$ are both true.\\ 
\textbf{\textit{Inductive step:}} Assume that $x_{1_n}^n(i-1),\ldots,x_{(i-1)_n}^n(i-1)$ are the variables verifying~\ref{h1}$(i-1)$ and \ref{h2}$(i-1)$ at iteration $i-1<K$. Let the variables at iteration $i$ be $x_{1_n}^n,\ldots,x_{i_n}^n$. We consider two cases:
\begin{enumerate}[leftmargin=*,label=\roman*)]
\item We first suppose that: 
\begin{equation}\label{ap:eq:SCPC_case1}
x^* = \textsc{Argmax}\tilde{f}(i,i, \mathcal{I}^n,\bar{P}^n) \leq x_{(i-1)_n}^n(i-1).    
\end{equation} 
In this case, Algorithm~\ref{algoSCPC} sets $x_{i_n}^n = x^*$ and $x_{l_n}^n = x_{l_n}^n(i-1)$, for all $l < i$. The induction hypothesis~\ref{h2}$(i-1)$ states that $\bar{P}^n \geq x_{1_n}^n \geq \cdots \geq x_{(i-1)_n}^n \geq 0$. By taking into account Eqn.~\eqref{ap:eq:SCPC_case1}, this inequality becomes $\bar{P}^n \geq x_{1_n}^n \geq \cdots \geq x_{(i-1)_n}^n \geq x^* = x_{i_n}^n \geq 0$.
Thus, \ref{h2}$(i)$ is satisfied. In addition, we know from~\ref{h1}$(i-1)$ that $x_{1_n}^n,\ldots,x_{(i-1)_n}^n$ maximizes $\sum_{l=1}^{i-1}{\tilde{f}_{l,l}^n}$. Since, the objective is separable and $x_{i_n}^n = x^*$ maximizes $\tilde{f}_{i,i}^n$ by construction, \ref{h1}$(i)$ is true.

\item Now, suppose that we have the opposite: 
\begin{equation}\label{ap:eq:SCPC_case2}
x^* = \textsc{Argmax}\tilde{f}(i,i, \mathcal{I}^n,\bar{P}^n) > x_{(i-1)_n}^n(i-1).    
\end{equation}  
In this case, the allocation mentioned above would violate constraint $C2'\text{--}3'$. 
The algorithm finds the highest index $j \in\{1,\ldots,i-2\}$ such that $x_{j_n}^n(i-1) \geq \textsc{Argmax}\tilde{f}(j+1,i,\mathcal{I}^n,\cU'_n,\bar{P}^n)$ in the while loop at line~6. Such an index exists since all variables are upper bounded by $\bar{P}^n$ and $x_{1_n}^n = \bar{P}^n$ due to Lemma~\ref{le:maxfij}. Let us show by contradiction that~\ref{h1}$(i)$ and~\ref{h2}$(i)$ are only satisfied if $x_{(j+1)_n}^n=\cdots=x_{i_n}^n$. If it is not the case, let $k > j+1$ be the last index such that $x_{k_n}^n=x_{(k+1)_n}^n= \cdots =x_{i_n}^n$ and $x_{(k-1)_n}^n>x_{k_n}^n$. We know from the while condition that $x_{(k-1)_n}^n<x^{*\prime}$, with $x^{*\prime}=\textsc{Argmax}\tilde{f}(k,i,\mathcal{I}^n,\cU'_n,\bar{P}^n)$. According to Lemma~\ref{le:maxfij}, $\tilde{f}_{k,i}^n$ is increasing on $\left[0,x^{*\prime}\right]$. Therefore, we can improve the objective function by setting $x_{k_n}^n, \ldots, x_{i_n}^n \gets x_{(k-1)_n}^n$. This is a contradiction with $x_{(k-1)_n}^n>x_{k_n}^n$, we have thus $x_{(j+1)_n}^n=\cdots=x_{i_n}^n$.
Furthermore, at the termination of the while loop, we have $\textsc{Argmax}\tilde{f}(j+1,i,\mathcal{I}^n,\cU'_n,\bar{P}^n) \leq x_{j_n}^n(i-1)$, which can be treated as in case~i). Hence, variables $x_{(j+1)_n}^n,\ldots,x_{i_n}^n$ are set equal to $\textsc{Argmax}\tilde{f}(j+1,i,\mathcal{I}^n,\cU'_n,\bar{P}^n)$ at line~10, and it satisfies~\ref{h1}$(i)$ and~\ref{h2}$(i)$.
\end{enumerate}

We proved that, in both cases~i) and~ii), the allocation $x_{1_n}^n,\ldots,x_{i_n}^n$ computed by Algorithm~\ref{algoSCPC} satisfies~\ref{h1}$(i)$ and \ref{h2}$(i)$. Therefore, by mathematical induction, the allocation returned at line~12 satisfies~\ref{h1}$(|\cU'_n|)$ and \ref{h2}$(|\cU'_n|)$. We note that~\ref{h1}$(|\cU'_n|)$ and \ref{h2}$(|\cU'_n|)$ are equivalent to an optimal solution of~\ref{PsubSCPC1}, which concludes the proof.
\null \hfill\ensuremath{\square}

\subsection{Proof of Theorem~\ref{th:algoiSCPC}}\label{ap:algoiSCPC}

\textbf{\textit{Optimality analysis:}} 
Let $x_{1_n}^n,\ldots,x_{|\cU'_n|_n}^n$ be the optimal allocation of \textsc{SCPC} with budget $P_{max}$. We consider now a lower budget $\bar{P}^n \leq P_{max}$. At each iteration $i$ of the loop in \textsc{SCPC}$\left(\mathcal{I}^n, \cU'_n, \bar{P}^n\right)$, the value $\textsc{Argmax}\tilde{f}\!\left(j,i,\mathcal{I}^n,\cU'_n,\bar{P}^n\right)$ can be replaced by $\min\{\textsc{Argmax}\tilde{f}\!\left(j,i,\mathcal{I}^n,\cU'_n,P_{max}\right),\bar{P}^n\}$, since they are equal by definition.
One can show, by mathematical induction on $i_n$, that the function \textsc{SCPC}$\left(\mathcal{I}^n, \cU'_n, \bar{P}^n\right)$ returns $\min\{x_{1_n}^n,\bar{P}^n\},\ldots,\min\{x_{|\cU'_n|_n}^n,\bar{P}^n\}$. Therefore, the latter allocation is also optimal.

\textbf{\textit{Complexity analysis:}} The initialization consists in running \textsc{SCPC}, with complexity $O\!\left(M^2\right)$ (see Theorem~\ref{th:algoSCPC}). Each subsequent evaluation requires to compute $\min\{x_{i_n}^n,\bar{P}^n\}$, for $i\in\{1,\ldots,|\cU'_n|\}$, with complexity $O\!\left(M\right)$.
\null \hfill\ensuremath{\square}

\subsection{Proof of Theorem~\ref{th:algoSCUS}}\label{ap:thalgoSCUS}

\textbf{\textit{Complexity analysis:}} The complexity mainly comes from the computation of $V$, $X$ and $U$ in the for loop from lines~13 to~27, which requires $M\sum_{i=1}^{K-1}\!\left(i\right) = O\!\left(MK^2\right)$ iterations. Each iteration has a constant number of operations. Thus, the overall worst case computational complexity is $O\!\left(MK^2\right)$.

\textbf{\textit{Optimality analysis:}} We will prove by induction that at any iteration $m\in\{0,\ldots,M\},\;j\in\{1,\ldots,K\}$ and $i\geq j$ of Algorithm~\ref{alg3}, the construction of $V[m,j,i]$ is the optimal value of problem~\ref{defV}. It follows directly that $V[M,1,1]$ is the optimal value of~\ref{Psub2}.\\
\textbf{\textit{Basis:}} For $m=0$, no user can be active due to constraint $C4'$. Thus, $V[0,j,i] = f_{j,K}^n(0)$ and $X[0,j,i]$ is initialized to zero. Furthermore, $U[0,j,i]=\varnothing$ to indicate that there is no previous index in the recursion. For simplicity of the algorithm, $V, X, U$ are also initialized for $j \leq i=K$ as explained in Section~\ref{sec:SCUS}.\\
\textbf{\textit{Inductive step:}} Let $m\in\{1,\ldots,M\}$ and $1\leq j \leq i\leq K-1$. Assume that $V[m',j',i']$ is the optimal value of~$\cP'_{SC}[m',j',i']$ for any $m'\leq m$, $j'\geq j$ and $i' > i$. We denote the optimal solution of problem~\ref{defV} by $x_j^n,\ldots,x_{K}^n$. 
Let $v_{act}$ (resp. $v_{inact}$) be the optimal value of~\ref{defV}, given that user $i$ is active (resp. inactive).
Let $x_{(i+1)_n}^{n*} = X[m-1,i+1,i+1]$ be the optimal value of $x_{(i+1)_n}^n$ in $\cP'_{SC}[m-1,i+1,i+1]$. If $x^* \leq x_{(i+1)_n}^{n*}$, then we can prove as in case~ii) of Appendix~\ref{ap:thalgoSCPC}, that user $i$ is inactive in the optimal solution. In this case, $V[m,j,i] = v_{inact}$. Otherwise, the optimal is $V[m,j,i] = \max\{v_{act},v_{inact}\}$.
Values $v_{act}$ and $v_{inact}$ are computed as follows:
\begin{itemize}[leftmargin=*]
\item Case $v_{inact}$: Suppose that the optimal solution of problem~\ref{defV} is achieved when user $i$ is inactive, then we have $x_{i}^n = x_{i+1}^n$ by definition of $\cU'_n$. It follows from $C5'$ that $x_j^n = \cdots = x_{i+1}^n$. We obtain, by definition, $V[m,j,i] = V[m,j,j+1]$, which we denote by $v_{inact}$.
\item Case $v_{act}$: Suppose now that user $i$ is active. Since $x^* > x_{(i+1)_n}^{n*}$ satisfies $C3'$, and the objective is separable, the optimal is obtained when maximizing independently $f_{j,i}^n$ and $\sum_{l=i+1}^{K}{f_l^n}$ with $m-1$ active users. That is, $V[m,j,i] = v_{act} \triangleq f_{j,i}^n\left(x^*\right) + V[m-1,i+1,i+1]$, where $x^*= \textsc{Argmax}f\!\left(j,i,\mathcal{I}^n,\bar{P}^n\right)$ in line~15.
\end{itemize}
Hence, $V[m,j,i]$, as computed in~\eqref{SCUSrec}, corresponds to the optimal of~\ref{defV}.

We derive, by mathematical induction, that $V[M,1,1]$ is the optimal value of $\cP'_{SC}[M,1,1]=$\ref{Psub2}. The corresponding optimal allocation $\bx^n$ is retrieved in lines~28 to~35.\hfill\ensuremath{\square}

\subsection{Proof of Theorem~\ref{th:algoiSCUS}}\label{ap:algoiSCUS}

\textbf{\textit{Optimality analysis:}} Let $y_1^{n},\ldots,y_K^{n}$ be the optimal solution of~\ref{Psub2} subject to a power constraint $\bar{P}^n$. Let $i\in\{1,\ldots,K\}$ be the unique index such that $y_1^n = \cdots = y_i^n$ and $y_i^n > y_{i+1}^n$. We know from Lemma~\ref{le:maxfij} that $y_1^n = \cdots = y_i^n = \bar{P}^n$. Therefore, $y_{i+1}^n,\ldots,y_K^n$ are all strictly less than $\bar{P}^n$. Let $x_1^{n},\ldots,x_K^{n}$ be the optimal solution of~$\cP'_{SC}[M,1,i]$ in the execution of \textsc{SCUS}$\left(\mathcal{I}^n, M, P_{max}\right)$, i.e., subject to a power budget $P_{max}$. According to Lemma~\ref{le:maxfij}, $x_1^n = \cdots = x_i^n = P_{max}$. We deduce from $f$'s unimodality in Lemma~\ref{le:maxfij}, that $y_{i+1}^n,\ldots,y_K^n$ is the optimal solution of~$\cP'_{SC}[M,i+1,i+1]$ given any power budget no less than $\bar{P}^n$. In particular, we have $x_l^n=y_l^n$, for all $l\in\{i+1,\ldots,K\}$. Hence, $x_1^{n},\ldots,x_K^{n}$ and $y_1^{n},\ldots,y_K^{n}$ correspond to the same user selection $\cU'_n$, and we derive $y_{l_n}^n = \min\{x_{l_n}^n,\bar{P}^n\}$, for $1 \leq l \leq |\cU'_n|$.

We proved above that, for any $\bar{P}^n \leq P_{max}$, there exists $\left(\cU'_n, x_{1}^n,\ldots,x_{K}^n\right)$ in $collection$, such that the optimal allocation subject to the power constraint $\bar{P}^n$ is $\min\{x_{1_n}^n,\bar{P}^n\},\ldots,\min\{x_{|\cU'_n|_n}^n,\bar{P}^n\}$.
Thus, the optimal user selection and power control is the one maximizing $F^n(\cU'_n,\bar{P}^n) = \sum_{l=1}^{|\cU'_n|} \tilde{f}_{l,l}^n\left(\cU'_n,\min\{x_{l_n}^n,\bar{P}^n\}\right) + B^n$ over all elements in $collection$, as shown at line~6.

\textbf{\textit{Complexity analysis:}} The initialization consists in running \textsc{SCUS}, with complexity $O\!\left(MK^2\right)$ (see Theorem~\ref{th:algoSCUS}). Each subsequent evaluation has complexity $O\!\left(MK\right)$. Indeed, there are $K$ active users sets $\cU'_n$ in $collection$, one for each solution of~$\cP'_{SC}[M,1,i]$, for $i\in\{1,\ldots,K\}$. For each of the $K$ possible active users set $\cU'_n$ in $collection$, we compute $F^n(\cU'_n,\bar{P}^n)$ with complexity $O\!\left(|\cU'_n|\right)=O\!\left(M\right)$.
\hfill\ensuremath{\square}

\subsection{Proofs of Lemma~\ref{le:derivative} and Theorem~\ref{th:algoGradJSPA}}\label{ap:derivative}

Let $x_{1}^n,\ldots,x_{K}^n$ be the output of \mbox{i-\textsc{SCUS}$\left(\bar{P}^n\right)$}, and $\cU'_n$ the corresponding active users set. For $i\in\{1,\ldots,K\}$, there exists $q \leq i$ and $q' \geq i$, such that $x_{q_n}^n=\cdots=x_{q'_n}^n=\textsc{Argmax}\tilde{f}\!\left(q,q',\mathcal{I}^n,\cU'_n,\bar{P}^n\right)$, according to Lemma~\ref{le:propOptimal_sumf}. We have:
\begin{align*}
\tilde{f}_{q,q'}^{n}\!\left(\cU'_n,\min\{x_{q_n}^n,\bar{P}^n\}\right) = \begin{cases}
\tilde{f}_{q,q'}^{n}\!\left(\cU'_n,\bar{P}^n\right), &\!\! \text{if $\bar{P}^n \leq x_{q_n}^n$},\\
\tilde{f}_{q,q'}^{n}\!\left(\cU'_n,x_{q_n}^n\right), &\!\! \text{if $\bar{P}^n > x_{q_n}^n$}.
\end{cases}
\end{align*}
We consider it as a function of $\bar{P}^n$. Its left derivative at $\bar{P}^n = x_{q_n}^n$ is $0$, according to Lemma~\ref{le:maxfij}. Its right derivative at $\bar{P}^n = x_{q_n}^n$ is $0$, as it is constant for $\bar{P}^n > x_{q_n}^n$. Hence, $\tilde{f}_{q,q'}^{n}\!\left(\cU'_n,\min\{x_{q_n},\cdot\}\right)$ is continuously differentiable on $\left[0,P_{max}\right]$.

Let $l$ be the greatest index such that $x_l^n = \bar{P}^n$. The function $F^n(\cU'_n,\bar{P}^n)$ can be written as $f_{1,l}^n\left(\bar{P}^n\right) + \sum_{i=l+1}^{K} f_{i,i}^n\left(x_{i}^n\right) + B^n$. Its derivative can be obtained by applying Eqn.~\eqref{proof:first_deriv_f} of Appendix~\ref{ap:propf} as follows:
\begin{equation}\label{ap:eq:derivative}
F^{n\prime}\!\left(\cU'_n,\bar{P}^n\right) = f_{1,l}^{n\prime}\!\left(\bar{P}^n\right) = \frac{W_n w_{\pi_n(l)}}{\left(\bar{P}^n+\tilde{\eta}^n_{\pi_n(l)}\right)\ln\!{(2)}}.
\end{equation}
As $F^n(\bar{P}^n) = \max_{\,\cU'_n}\{F^n(\cU'_n,\bar{P}^n)\}$, where $\max$ is taken over all active users sets in $collection$ of \mbox{i-\textsc{SCUS}}, and the $\max$ operator only preserves semi-differentiability, Eqn.~\eqref{ap:eq:derivative} is the left derivative of $F^{n}$. This proves Lemma~\ref{le:derivative}. 

In addition, the second left derivative of $F^{n}$ satisfies:
\begin{equation}\label{ap:eq:derivative2}
\beta
\leq F^{n\prime\prime}\!\left(\bar{P}^n\right) 
= \frac{-W_n w_{\pi_n(l)}}{\left(\bar{P}^n+\tilde{\eta}^n_{\pi_n(l)}\right)^2\ln\!{(2)}}
\leq \alpha < 0,
\end{equation}
where $\beta$ and $\alpha$ are constant and defined as:
\begin{equation*}
\beta = \frac{-W_n w_{\pi_n(l)}}{\left(\tilde{\eta}^n_{\pi_n(l)}\right)^2\ln\!{(2)}}, \, 
\alpha = \frac{-W_n w_{\pi_n(l)}}{\left(P_{max}+\tilde{\eta}^n_{\pi_n(l)}\right)^2\ln\!{(2)}}.
\end{equation*}
Although $F^n$ is only semi-differentiable at some points, it is twice differentiable on each interval where the optimal user selection $\cU'_n$ does not change. Appendix~\ref{ap:algoiSCUS} shows that there are $K$ such intervals. Eqn.~\eqref{ap:eq:derivative2} implies that $F^{n}$ is piece-wise twice differentiable, $\alpha$-strongly concave and $\beta$-smooth. Therefore, the projected gradient descent on the simplex $\mathcal{F}$ converges in $O\!\left(\log\!\left(1/\xi\right)\right)$ iterations, according to~\cite[Section 2.2.4]{nesterov1998introductory}.
This proves Theorem~\ref{th:algoGradJSPA}. \hfill\ensuremath{\square}

\subsection{Proof of Theorem~\ref{th:gap}}\label{ap:gap}
Let $\bar{P}^{n*}$ be the optimal power budget of~\ref{Psub1} on subcarrier $n\in\cN$. The power budget after discretization with step size $\delta$ is denoted by $a_n^* \triangleq \lfloor\bar{P}^{n*}/\delta\rfloor\delta$. We have:
\begin{align}
F^* - &F^*_{MCKP} \leq \sum_{n\in\cN}\left(F^n\!\left(\bar{P}^{n*}\right) - F^n\!\left(a_n^*\right)\right), \label{ap:gap_eq1}\\
&\leq \sum_{n\in\cN}\max_{\,\cU'_n}\{F^{n\prime}(\cU'_n,\bar{P}^{n*})\}\times\left(\bar{P}^{n*}-a_n^*\right), \label{ap:gap_eq2}\\
&\leq \delta\sum_{n\in\cN} \max_{k\in\cK}\Bigg\{\frac{W_n w_{\pi_n(k)}}{\left(\bar{P}^{n*}+\tilde{\eta}^n_{\pi_n(k)}\right)\ln\!{(2)}}\Bigg\}. \label{ap:gap_eq3}
\end{align}
Inequality~\eqref{ap:gap_eq1} comes from the definition of $F^*$ and the fact that $F^*_{MCKP} \geq \sum_{n\in\cN}F^n\!\left(a_n^*\right)$, as $F^*_{MCKP}$ is the optimal discrete solution with step size $\delta$.
We know from Appendix~\ref{ap:derivative} that $F^n(\bar{P}^n) = \max_{\,\cU'_n}\{F^n(\cU'_n,\bar{P}^n)\}$, and that $F^{n}\!\left(\cU'_n,\bar{P}^n\right)$ is twice differentiable and concave, for any $\cU'_n$ and $n\in\cN$. Hence, $F^n$ lies below the maximum slope tangent among the tangents of $F^n(\cU'_n,\bar{P}^{n*})$, for all $\cU'_n$. This implies inequality~\eqref{ap:gap_eq2}. We obtain~\eqref{ap:gap_eq3} by applying Eqn.~\eqref{ap:eq:derivative}, and the fact that $\bar{P}^{n*}-a_n^* \leq \delta$ by construction.
\hfill\ensuremath{\square}

\subsection{Proof of Theorem~\ref{th:algoOptJSPA}}\label{ap:algoOptJSPA}
Let us first briefly explain the principle of dynamic programming by weights. Let $Z$ be a 2D-array such that $Z[n,l]$ is defined as the optimal value of~\ref{P_MCKP} restricted to the first $n$ classes and with restricted capacity $l\cdot\delta$. It is initialized as $Z[0,l] = 0$, for any $l=0,\ldots,J$. For $n\in\cN$ and $l=0,\ldots,J$, the recurrence relation is:
\begin{equation*}
Z[n,l] = \max_{l'\leq l}\{Z[n-1,l-l']+c_{n,l'}\}.
\end{equation*}
The complexity and optimality proofs of \textsc{Opt-JSPA} are presented below.

\textbf{\textit{Optimality analysis:}} Reference~\cite{book:knapsack} proves that dynamic programming by weights is optimal for~\ref{P_MCKP}. Since problems~\ref{Psub1} and~\ref{P_MCKP} are equivalent, the proposed \textsc{Opt-JSPA} based on dynamic programming by weights is also optimal for~\ref{Psub1}.

\textbf{\textit{Complexity analysis:}} In Algorithm~\ref{algoOptJSPA}, we first transform~\ref{Psub1} to problem~\ref{P_MCKP}: from line~1 to~5, every item's profit $c_{n,l}$ is computed using i-\textsc{SCUS} in $O\!\left(NMK^2+JNMK\right)$. Then, we perform dynamic programming by weights at lines~6-7. According to~\cite{book:knapsack}, its complexity is $O\!\left(J^2N\right)$, which is the number of items $N\left(J+1\right)$ multiplied by the number of possible power values $J+1$. Therefore, the overall complexity is $O\!\left(NMK^2+JNMK+J^2N\right)$.
\hfill\IEEEQEDhere

\subsection{Estimation $U$ in Algorithm~\ref{algoepsilonJSPA}}\label{ap:estimationU}

In this section, we denote by $F^*_{MCKP}(P_{max})$ the optimal value of~\ref{P_MCKP} with cellular power budget $P_{max}$. We provide some properties in Lemma~\ref{le:propFopt} that will be used for the analysis of the estimation procedure.
\begin{myle}[Monotonicity and sublinearity of  $F^*_{MCKP}$]\label{le:propFopt}
$ $\newline
$F^*_{MCKP}$ is a non-decreasing and sublinear function of $P_{max}$. That is, for any $P_1 < P_2$, $F^*_{MCKP}(P_1) \leq F^*_{MCKP}(P_2)$ and $F^*_{MCKP}(P_1+P_2) \leq F^*_{MCKP}(P_1) + F^*_{MCKP}(P_2)$.
\end{myle}
\begin{IEEEproof}
We first prove the monotonicity of $F^*_{MCKP}$. Let $\mathcal{F}_1'$ and $\mathcal{F}_2'$ be two feasible sets of~\ref{Psub1} with power budget $P_1$ and $P_2$ respectively. Assuming $P_1 < P_2$, then any solution of $\mathcal{F}_1'$ is also a solution of $\mathcal{F}_2'$, i.e., $\mathcal{F}_1' \subset \mathcal{F}_2'$. Since~\ref{Psub1} is a maximization problem over $\mathcal{F'}$, we have $F^*_{MCKP}(P_1) \leq F^*_{MCKP}(P_2)$. This proves that $F^*_{MCKP}$ is non-decreasing.

Now, let us tackle the sublinearity of $F^*_{MCKP}$. 
We first prove that the $f_{j,i}^n$ are sublinears. If $j = 1$ or $w_{\pi_n(i)} \geq w_{\pi_n(j-1)}$, then $f_{j,i}^n$ is concave according to Lemma~\ref{le:maxfij}. Therefore, it is also sublinear. Otherwise, $f_{j,i}^n$ is concave before $c_2$ and decreasing after $c_1 \leq c_2$. In this case, $f_{j,i}^n$ is thus also sublinear. Secondly, for any subcarrier $n$ and user selection $\cU'_n$, \ref{PsubSCPC1} consists in maximizing a sum of separable sublinear functions $f^n_{j,i}$ subject to a budget constraint $\bar{P}^n$. Hence, $F^n\!\left(\cU'_n,\bar{P}^n\right)$ is sublinear in $\bar{P}^n$. Thirdly, the optimal of~\ref{Psub2} can be seen as the best allocation over all possible user selections, i.e., $F^n(\bar{P}^n) = \max_{\,\cU'_n}\{F^n\!\left(\cU'_n,\bar{P}^n\right)\}$. The $\max$ operator preserves sublinearity. Therefore, $F^n(\bar{P}^n)$ is sublinear in $\bar{P}^n$. Finally, $F^*_{MCKP}$ is sublinear in $P_{max}$, since~\ref{Psub1} is a separable sum maximization of $F^n$ subject to budget constraint $P_{max}$.
\end{IEEEproof}

Let us introduce a variant of~\ref{P_MCKP}, denoted by \textsc{MCKP}'. The differences are as follows. Its cellular power budget is $2P_{max}$. The item's weights can only take value of the form $a_{n,l} = l\lfloor J/N\rfloor\delta$ for $n\in\cN$, $l\in\{0,\ldots,2N\}$. The profits values are defined similarly as $c_{n,l} = F^n\!\left(a_{n,l}\right)$. Consequently, \textsc{MCKP}' only contains $2N+1$ items per class.
The idea of the proof is to show that a greedy solution of \textsc{MCKP}' is a constant factor approximation of~\ref{P_MCKP} optimal value. The value of $U$ is then easily obtained using the greedy Dyer-Zemel algorithm~\cite[Section 11.2]{book:knapsack}. In this case, the complexity is independent of $J$ and negligible compared to the rest of the algorithm. One could also get an estimation by applying the Dyer-Zemel algorithm directly to~\ref{P_MCKP}. However, the complexity would be proportional to $O\!\left(J\right)$ which is against the idea of polynomial-time approximation.

Let $y^*_{n,l}$, for $n\in\cN$, $l\in\{0,\ldots,2N\}$, be an optimal solution of this problem. In addition, we denote by $y'_{n,l}$ for $n\in\cN$, $l\in\{0,\ldots,2N\}$, a $1\!/2$-approximation given by the Dyer-Zemel algorithm.
On the one hand, we have:
\begin{align}
\sum_{n\in\cN}&\sum_{l=1}^{2N}{c_{n,l}y'_{n,l}} \geq \frac{1}{2}\sum_{n\in\cN}\sum_{l=1}^{2N}{c_{n,l}y^*_{n,l}}, \label{ap::estimationU_eqA}\\
&\geq \frac{1}{2}\sum_{n\in\cN}F^n\!\left(\Biggl\lceil\frac{\bar{P}^{n*}}{\lfloor J/N \rfloor\delta}\Biggr\rceil
\Biggl\lfloor\frac{J}{N}\Biggr\rfloor\delta
\right), \label{ap::estimationU_eqB}\\
&\geq \frac{1}{2}\sum_{n\in\cN}F^n\!\left(\bar{P}^{n*}\right) = \frac{1}{2}F^*_{MCKP}\!\left(P_{max}\right), \label{ap::estimationU_eqC}
\end{align}
where $\bar{P}^{n*}$ is the power allocated to subcarrier $n$ in $F^*_{MCKP}\!\left(P_{max}\right)$.
The $1\!/2$-approximation of $y'_{n,l}$ translates into Eqn.~\eqref{ap::estimationU_eqA}. The right term of Eqn.~\eqref{ap::estimationU_eqB} corresponds to a valid allocation of \textsc{MCKP}', with item $l = \lceil \bar{P}^{n*}/ (\lfloor J/N \rfloor\delta)\rceil$ allocated in class $n$. Indeed, by definition of the ceiling and floor functions, we have:
\begin{equation*}
\bar{P}^{n*} \stackrel{\text{(e)}}{\leq} \Biggl\lceil\frac{\bar{P}^{n*}}{\lfloor J/N \rfloor\delta}\Biggr\rceil
\Biggl\lfloor\frac{J}{N}\Biggr\rfloor\delta < \bar{P}^{n*}+ \Biggl\lfloor\frac{J}{N}\Biggr\rfloor\delta
\leq \bar{P}^{n*}+\frac{P_{max}}{N}.
\end{equation*}
Therefore, 
\begin{equation*}
\sum_{n\in\cN}{\Biggl\lceil\frac{\bar{P}^{n*}}{\lfloor J/N \rfloor\delta}\Biggr\rceil
\Biggl\lfloor\frac{J}{N}\Biggr\rfloor\delta} < \sum_{n\in\cN}\!\left(\bar{P}^{n*}+\frac{P_{max}}{N}\right) = 2P_{max}. 
\end{equation*}
In other words, the power budget is also satisfied. As it is a valid allocation for \textsc{MCKP}', it must have a total profit not greater than the optimal profit $\sum_{n\in\cN}\sum_{l=1}^{2N}{c_{n,l}y^*_{n,l}}$, which proves inequality~\eqref{ap::estimationU_eqB}. We derive Eqn.~\eqref{ap::estimationU_eqC} from inequality~(e) and the monotonicity of $F^n$  (see Lemma~\ref{le:propFopt}).

We have, on the other hand:
\begin{align}
\sum_{n\in\cN}\sum_{l=1}^{2N}{c_{n,l}y'_{n,l}} &\leq \sum_{n\in\cN}\sum_{l=1}^{2N}{c_{n,l}y^*_{n,l}}, \label{ap::estimationU_eq1}\\
&\leq F^*_{MCKP}\!\left(2P_{max}\right), \label{ap::estimationU_eq2}\\
&\leq 2F^*_{MCKP}\!\left(P_{max}\right). \label{ap::estimationU_eq3}
\end{align}
The optimality of $y^*_{n,l}$ implies Eqn.~\eqref{ap::estimationU_eq1}. Eqn.~\eqref{ap::estimationU_eq2} comes from the fact that the items of \textsc{MCKP}' is a subset of~\ref{P_MCKP} items, given a budget $2P_{max}$. Eqn.~\eqref{ap::estimationU_eq3} follows from the sublinearity of $F^*_{MCKP}$ (see Lemma~\ref{le:propFopt}).

Let $U \triangleq 2\sum_{n\in\cN}\sum_{l=1}^{2N}{c_{n,l}y'_{n,l}}$. We derive from inequalities~\eqref{ap::estimationU_eqC} and~\eqref{ap::estimationU_eq3} the desired approximation bound: 
\begin{equation*}
U \geq F^*_{MCKP}\!\left(P_{max}\right) \geq U/4. \IEEEQEDhereeqn
\end{equation*}

\subsection{Proof of Theorem~\ref{th:algoEpsJSPA}}\label{ap:algoEpsJSPA}

\textbf{\textit{Complexity analysis:}} 
We divide the complexity analysis of Algorithm~\ref{algoepsilonJSPA} in four parts as follows. The overall complexity can be obtained by summing the complexity of each part.\\
\textbf{\textit{i. Precomputation:}} The precomputation required for setting up i-\textsc{SCUS} on each subcarrier has complexity $O\!\left(NMK^2\right)$.\\
\textbf{\textit{ii. Line 1:}} 
The estimation procedure presented in Appendix~\ref{ap:estimationU}, consists in $O\!\left(N^2\right)$ function evaluations and $O\!\left(N^2\right)$ iterations of the Dyer-Zemel algorithm. Each function evaluation is computed by i-\textsc{SCUS}, therefore the complexity of this part is $O\!\left(N^2MK\right)$.\\
\textbf{\textit{iii. Lines 2-4:}} 
Each $L_n$, for $n\in\cN$, is obtained by multi-key binary search~\cite{tarek2008multi}. For each $L_n$, we need to find $4N/\text{\textepsilon}$ keys in an array $\{c_{n,1},\ldots,c_{n,J}\}$ of length $J$. Since repetition is not allowed, the binary search returns at most $\min\{4N/\text{\textepsilon},J\}$ items. More precisely, it computes each of the $4N/\text{\textepsilon}$ keys in time $\log(J)$, with at most $J$ function evaluations in total. Therefore, the binary search performs $O\!\left(\min\{\log\!\left(J\right)N/\text{\textepsilon},J\}\right)$ function evaluations. Multiplied by the complexity of each function evaluation on each subcarrier, we obtain $O\!\left(\min\{\log\!\left(J\right)N^2MK/{\text{\textepsilon}},JNMK\}\right)$.\\
\textbf{\textit{iv. Lines 5-6:}} 
Let us first briefly explain the dynamic programming by profits~\cite{book:knapsack}. Let $Y$ be the DP array such that $Y[n,q]$ denotes the minimal weight, i.e., minimal power budget, required to achieve WSR $q\cdot\text{\textepsilon}U/4N$ when problem~\ref{P_MCKP} is restricted to the first $n$ classes. It is initialized as $Y[0,0] = 0$ and $Y[0,q] = +\infty$, for $q=1,\ldots,\lfloor 4N/\text{\textepsilon} \rfloor$. For $n\in\cN$ and $q=0,\ldots,\lfloor 4N/\text{\textepsilon} \rfloor$, the recurrence relation is:
\begin{equation}\label{rec:dp_by_profits}
Y[n,q] = \min_{l\in L_n}\begin{cases}
Y\left[n-1,q-\Bigl\lfloor\frac{4c_{n,l}N}{\text{\textepsilon}U} \Bigr\rfloor\right]+a_{n,l}, \\
\hphantom{Y\left[n,q-\Bigl\lfloor\frac{c_{n,l}N}{\text{\textepsilon}U} \Bigr\rfloor\right]}\text{if $\frac{q\cdot\text{\textepsilon}U}{4N} \geq c_{n,l}$}, \\
+\infty, \text{ otherwise}.
\end{cases}
\end{equation}
This recursion has complexity $O\!\left(\min\{N^3/\text{\textepsilon}^2,J^2N\}\right)$, which is the number of all considered items $\sum_{n\in\cN}{|L_n|} = \min\{4N^2/\text{\textepsilon},JN\}$ multiplied by the number of comparisons in Eqn.~(\ref{rec:dp_by_profits}), $|L_n|=\min\{4N/\text{\textepsilon},J\}$. 

\textbf{\textit{Approximation analysis:}} As proved in~\cite[Section 11.9]{book:knapsack}, the optimal solution obtained by dynamic programming by profits considering only items in $L_n$, differs from $F^*_{MCKP}$ by at most a factor $1-\text{\textepsilon}$.

In summary, \textsc{\textepsilon-JSPA} achieves \textepsilon-approximation with polynomial complexity in $1/\text{\textepsilon}$ and $N$, $M$, $K$. Therefore, \textsc{\textepsilon-JSPA} is a FPTAS, which concludes the proof.
\hfill\IEEEQEDhere

\end{document}